\documentclass[12pt,letterpaper]{article}
\usepackage{renstyles}
\hypersetup{
	citecolor=green!60!black,  %
}

\usepackage{cases}

\usepackage{tikz,pgfplots}
\usetikzlibrary{pgfplots.groupplots}

\newcommand{\data}{Data} %

\definecolor{myblue}{RGB}{0, 114, 189}
\definecolor{myorange}{RGB}{216, 83, 25}
\definecolor{myyellow}{RGB}{237, 177,  32}
\definecolor{mypink}{RGB}{239, 71,  111}
\definecolor{mypurple}{RGB}{126,   47,  142}
\definecolor{mygreen}{RGB}{119,  172,  48}
\definecolor{mycharcoal}{RGB}{38, 70, 83}
\definecolor{mycyan}{RGB}{77,  190,  238}
\definecolor{mylightred}{RGB}{255,  204,  203}
\definecolor{mygray}{RGB}{128, 128, 128}
\newcommand{\mylinewidth}{1pt}

\usepackage{cite}

\pgfplotscreateplotcyclelist{cyclelist-compare}{
	blue!80!green, every mark/.append style={solid,fill=blue!80!white},mark=square*,line width = \mylinewidth\\
	red!90!black,every mark/.append
	style={solid,fill=red!80!white},mark=triangle*,line width =\mylinewidth\\ 
	myblue,,densely dashed, every mark/.append style={solid,fill=myblue!80!white},mark=*,line width = \mylinewidth\\
	myorange, densely dashed,every mark/.append style={solid,fill=myorange!70!white},mark=diamond*,line width =\mylinewidth\\
}

\pgfplotscreateplotcyclelist{cyclelist-reconstruct}{
	myyellow, line width = 1.5pt\\
	blue!80!green,every mark/.append
	style={solid,fill=red!80!white},mark=none, dashed, line width =1.2\\ 
	red!90!black, densely dotted, every mark/.append style={solid,fill=myblue!80!white},mark=none,line width = 1.2\\
	mygreen, densely dashed,every mark/.append style={solid,fill=myorange!70!white},mark=none,line width =1.2\\
}

\pgfplotscreateplotcyclelist{cyclelist-reconstruct-longtime}{
	myyellow, line width = 1.5pt\\
	blue!80!green,every mark/.append
	style={solid,fill=red!80!white},mark=none, dashed, line width =1.2\\ 
	myblue!80!cyan,every mark/.append
	style={solid,fill=red!80!white},mark=none, densely dashed, line width =1.2\\  
	red!70!black, dotted, every mark/.append style={solid,fill=myblue!80!white},mark=none,line width = 1.2\\
	red!70!white, densely  dotted, every mark/.append style={solid,fill=myblue!80!white},mark=none,line width = 1.2\\
	mygreen, densely dashed,every mark/.append style={solid,fill=myorange!70!white},mark=none,line width =1.2\\
}

\numberwithin{equation}{section}

\renewcommand{\bp}{p}

\renewcommand{\bx}{x} 
\renewcommand{\by}{y} 

\newcommand{\Du}{\nabla u}

\newcommand{\Dhu}{\nabla \hat{u}}
\newcommand{\Dhm}{\nabla \hat{m}}

\newcommand{\Dv}{\nabla v}

\newcommand{\DpH}{H_{\bp}}

\newcommand{\barU}{\bar{U}}
\newcommand{\barM}{\bar{M}}
\newcommand{\DbU}{\nabla \barU}
\newcommand{\DbM}{\nabla \barM}
\newcommand{\Db}{\nabla b}

\newcommand{\XTM}{X^T_{M_1}}
\newcommand{\cTK}{\cT_K}

\newcommand{\DataFun}{g}
\newcommand{\MeasOpr}{\cG}

\newcommand{\edit}[1]{#1}

\title{A Policy Iteration Method for Inverse Mean Field Games\thanks{Dedicated to Professor Bj\"{o}rn Engquist on the occasion of his 80th birthday!}}
\author{
	Kui Ren\thanks{Department of Applied Physics and Applied Mathematics, Columbia University, New York, NY 10027; kr2002@columbia.edu}
	\and
	Nathan Soedjak\thanks{Department of Applied Physics and Applied Mathematics, Columbia University, New York, NY 10027; ns3572@columbia.edu}
	\and
	Shanyin Tong\thanks{Department of Mathematics, University of Pennsylvania, Philadelphia, PA 19104; tong3@sas.upenn.edu}
}

\begin{document}

\maketitle

\begin{abstract}
    We propose a policy iteration method to solve an inverse problem for a mean-field game (MFG) model, specifically to reconstruct the obstacle function in the game from the partial observation data of value functions, which represent the optimal costs for agents. The proposed approach decouples this complex inverse problem, which is an optimization problem constrained by a coupled nonlinear forward and backward PDE system in the MFG, into several iterations of solving linear PDEs and linear inverse problems. This method can also be viewed as a fixed-point iteration that simultaneously solves the MFG system and inversion. We prove its linear rate of convergence. In addition, numerical examples in 1D and 2D, along with performance comparisons to a direct least-squares method, demonstrate the superior efficiency and accuracy of the proposed method for solving inverse MFGs.
\end{abstract}



\section{Introduction}
\label{SEC:Intro}

Mean-field games (MFGs) model the strategic interactions among a large number of rational agents, characterizing the Nash equilibrium in such games \cite{LaLi-JJM07}.
 They are particularly useful in situations when the collective behavior of all agents significantly influences outcomes while the influence of a single agent is negligible. 
Applications include, but are not limited to, modeling market behavior in economics and finance~\cite{AcHaLaLiMo-RES22, CaDeLa-AMO07}, capturing crowd dynamics like traffic flow~\cite{HuChDiDu-TRRC21}, and studying disease spread in epidemiology \cite{LeLiTeLiOs-SIAM21}. In addition, MFGs have recently gained significant attention because of their connections to machine learning~\cite{RuOsLiNuFu-PNAS20}, including reinforcement learning for dynamic optimization to obtain their solutions and deep learning to address the curse of dimensionality. Numerous studies also focus on numerical schemes for solving MFGs~\cite{AcCa-SIAM10,AgLeFuNu-JCP22,FuLiOsLi-arXiv23,LiJaLiNuOs-SIAM21,LeLiTeLiOs-SIAM21,YuXiChLa-arXiv24}, given their unique coupled structure of forward and backward nonlinear PDEs.

In MFGs, each agent selects a strategy to influence its state dynamics, aiming to minimize a cost function that typically depends on the agent's state, control actions, and the mean-field term (such as population density). The optimal control and corresponding distribution are the deterministic solutions of the MFG PDE system: a backward Hamilton-Jacobi-Bellman (HJB) equation for the value function (optimal cost of a single agent) and a Fokker-Planck (FP) equation for the population density.
Since the cost function may involve environmental information that is independent of the game and not always known from first principles or direct measurements, this leads to the inverse mean-field game problem, where the goal is to determine parts of the cost function or environmental information based on observations of MFG solutions, such as optimal costs, strategies, or densities.
This information helps us calibrate the models with data, enabling more accurate models in the future.

MFGs and their associated inverse problems are useful in practice and mathematically intriguing. However, their numerical and theoretical analysis face significant challenges. The primary difficulty arises from the coupled structure of the forward and backward nonlinear differential equations, which cannot be addressed using standard numerical methods that integrate time in a single direction. This complexity is further compounded in the context of inverse problems, as the coupled systems act as equality constraints on the optimization process, resulting in adjoint equations with a similar structure. Recently, there has been progress in this field, with research focusing on theoretical analysis \cite{DiLiZh-arXiv23,ImLiYa-IPI24,ImYa-arXiv23,KlAv-SIAM24,KlLiLi-JIIP24,LiMoZh-IP23, LiZh-arXiv24, ReSoWaZh-IP24}
and numerical schemes~\cite{DiLiOsYi-JSC22,ChFuLiNuOs-IP23,YuXiChLa-arXiv24,GuMoYaZh-JCP24} within specific settings. Despite these advancements, numerous open questions and unexplored scenarios remain.

In this paper, we aim to introduce a computationally efficient and theoretically guaranteed method for solving the inverse MFG problem. We focus on reconstructing the obstacle function $b(\bx)$~\cite{YuXiChLa-arXiv24}, %
reflecting environmental information and sometimes called the potential, %
from observations of the value function, assuming the system follows the MFG model \eqref{EQ:MFG}. For example, the observations could be the optimal cost values of an agent at the initial time.

To achieve our goals, we propose a policy iteration method for inverse MFG problems, inspired by its use in optimal control and MFG systems. This method decouples the MFG using an auxiliary function, the policy, representing the optimal strategy for minimizing the total cost. It simplifies the complicated coupled-nonlinear-PDE-constrained optimization into iterative steps:
(i) solve a linear Fokker-Planck equation;
(ii) solve a linear inverse problem constrained by a linearized HJB equation from the Legendre transform;
(iii) perform a pointwise optimization to update the policy, often in closed form. The convergence of these iterations is demonstrated theoretically and numerically. The decoupling also eliminates nonlinearity and nonconvexity, ensuring the efficiency of the policy iteration method for inverse MFGs, as demonstrated in comparisons with a standard least-squares PDE-constrained optimization method.

The main \emph{contributions} of this paper are: (1) we propose an efficient approach: policy iteration for inverse MFGs, which decouples the MFG system and only requires several iterations of solving linear PDEs and linear inverse problems; (2) we prove that our iterative approach achieves a linear rate of convergence; (3) we demonstrate the superior efficiency of this method compared to direct least-squares PDE-constrained optimization approaches for inverse problems, which require solving additional coupled PDE systems for gradient evaluation. This efficiency is illustrated through both one-dimensional and two-dimensional experiments.

\section{Preliminaries\label{sec:prelim}}

\subsection{Mean field games and their inverse problems\label{sec:MFG}}

In this paper, we study mean field games, which typically have a strongly coupled structure, consisting of a forward nonlinear Fokker-Planck (FP) equation capturing the evolution of the distribution of all agents in the game, and a nonlinear backward Hamilton-Jacobi-Bellman (HJB) equation characterizing how the total cost of an individual agent changes:
\begin{equation}\label{EQ:MFG}
\left\{\begin{array}{ll}
-\partial_t u-\eps\Delta u + H(\Du) = b(\bx) + F(m)\, &\text{in } Q,\\
\partial_t m - \eps\Delta m - {\rm div}(m \DpH(\Du)) = 0\, &\text{in }Q,\\
u(\bx,T) = u_T(\bx),\ m(\bx,0)=m_0(\bx)\, &\text{in } \bbT^d,
\end{array}\right.
\end{equation} 
where $Q:=\bbT^d \times [0,T]$, and the state space is a flat torus $\bbT^d:=\bbR^d/\bbZ^d$, which can be viewed as a unit cube $[0,1]^d$ with periodic boundary conditions.
$\bx\in\bbT^d$ denotes the state variable (for example, location in space), $T$ is the final time, $t\in[0,T]$ denotes the time variable. 
The notation $\Du$ and $\Delta u$ are, respectively, the gradient and Laplacian of $u$ with respect to $\bx$-variable, and ${\rm div}$ is the divergence operator with respect to $\bx$ as well. The diffusion parameter $\eps>0$ captures the intensity of random fluctuations. 
$u(\bx, t)\in \bbR$ is the value function representing the minimal cost of a single agent starting at time $t$ at state $\bx$; its evolution follows a backward HJB equation.
$m(\bx, t)\in \bbR_{\ge 0}$ stands for the density of all agents at state $\bx$ at time $t$; its evolution is governed by a forward FP equation. 
The initial condition $m_0:\bbT^d \to \bbR_{\ge 0}$ represents the initial distribution satisfying $\int_{\bbT^d} m_0(\bx) d\bx = 1$. 
The obstacle function (sometimes called potential) is $b: \bbT^d \to \bbR$, which only depends on the current state. The name ``obstacle" comes from the fact that agents are penalized for moving in regions where the value of $b$ is high.
Besides the obstacle function $b$, the total cost is also determined by the terminal cost function $u_T:\bbT^d \to \bbR$, the interaction cost function $F:\bbR \to \bbR$ (applying to the density, a usual choice is $F(m) =m^2$), and the dynamic cost captured by the Hamiltonian $H:\bbR^d \to \bbR$ (applying to $\bp=\Du$, with its gradient denoted as $\DpH$). 

Note that in a more general setting, the Hamiltonian can depend on $\Du$, $\bx$, and $m$ together to capture the total cost itself. In that case, there is no term on the right-hand side of the HJB equation in \eqref{EQ:MFG}. In our paper, we consider a separable setting where the costs can be split into parts: $H(\Du)$ for dynamics, $b(\bx)$ for obstacles, and $F(m)$ for interactions. This is inspired by the common choice for $H$ to represent the kinetic energy, 
$H(\bp) = \frac{1}{2}|\bp|^2$, where $|\cdot|$ denotes the norm of a vector in $\bbR^d$. 
Our work can be directly extended to problems with $H$ and $F$ also having a dependence on $\bx$. For clarity, we only discuss the case where they do not involve $\bx$.

The inverse MFG problem aims to reconstruct certain functions in~\eqref{EQ:MFG} using observational data of its solutions. While extensive research exists on solving MFG systems, studies on inverse MFGs are relatively limited. On the theoretical side, there has been considerable interest recently in proving injectivity results for forward maps in various setups and stability estimates for the inverse maps~\cite{DiLiZh-arXiv23,ImLiYa-IPI24,ImYa-arXiv23,KlAv-SIAM24,KlLiLi-JIIP24,LiMoZh-IP23, LiZh-arXiv24, ReSoWaZh-IP24}. The techniques used typically include either the linearization method of Isakov or the Carleman estimates. On the computational side, different numerical algorithms have been developed for computational inversion. 
A primal-dual method was introduced for solving inverse problems in mean-field games by formulating the data misfit minimization as a standard PDE-constrained optimization problem~\cite{DiLiOsYi-JSC22}, aiming to reconstruct metrics and interaction kernels from observations of density and velocity fields.
\cite{ChFuLiNuOs-IP23} developed an operator splitting algorithm to learn running costs and interaction energy from boundary observations of population density and strategies. \cite{YuXiChLa-arXiv24} proposed a bilevel optimization to reconstruct obstacles and metrics from distribution and strategy data. \cite{GuMoYaZh-JCP24} employed Gaussian processes to recover strategies and environmental configurations from population and partial environmental observations, and \edit{\cite{KlLiYa-AMO24}} presented a method based on Carleman estimates to convexify a mean field game inverse problem. \edit{In \cite{KlLiYa-CMA25} the convexification
method is applied to solve a coefficient inverse problem of MFGs, and
\cite{KlLi-DG25} summarizes results of its authors on applications of Carleman estimates to both
forward and inverse problems for the MFGs.}

In this work, we focus on developing an efficient algorithm that provably converges to a (not necessarily unique) solution to the inverse problem. We consider the following inverse MFG problem setup: Given that the system follows the MFG \eqref{EQ:MFG} and the functions $H, F, u_T$, and $m_0$, we seek to reconstruct the obstacle function $b(\bx)$ from data $\DataFun(\bx)$. Specifically, we assume the data is a linear measurement $\MeasOpr u$ of the value function, i.e., $\DataFun\approx \MeasOpr u$. In noiseless cases, $\DataFun=\MeasOpr u$.
In our numerical experiments, we consider two types of data: (i) the initial value $\MeasOpr u:=u(\bx, 0)$ and (ii) the derivative at the final time $\MeasOpr u:=\partial_t u(\bx, T)$.
\edit{Note that case (ii), which involves time-derivative data, is not the standard setting in inverse problems; however, such data can be obtained by finite-difference approximations of the value function over a short time interval near the terminal time.}

\subsection{Policy iteration for solving MFGs\label{sec:Policy-MFG}}

Policy iteration is a classical algorithm for optimal control problems; see~\cite{Bellman-Book57,Howard-Book60,Bertsekas-Book87,Bertsekas-JCTA11,FaGhMaSz-ANIPS08,LaPa-JMLR03,PuBr-MOR79,Scherrer-PMLR14,SaRu-SIAM04,SuBa-Book18} and references therein for some random samples of earlier works and recent developments in the subject. A policy iteration method for solving the MFG system~\eqref{EQ:MFG} was first introduced in \cite{CaCaGo-COCV21}, and studied further in \cite{CaTa-JMAA22, LaSoTa-AMO23, TaSo-SIAM24, AsMi-arXiv23}. 
It decouples the MFG using an auxiliary function called the policy $q$, 
which is the optimal strategy to the control problem minimizing the total cost (whose optimal value forms the value function $u$).
Let $L: \bbR^d\to \bbR$ be the Lagrangian (Legendre transform  of the Hamiltonian $H$), i.e., $L(q) := \sup_{p\in \bbR^d} (q\cdot p - H(p))$. Fix $R>0$ and choose a bounded vector field $q^{(0)}: Q\to \bbR^d$ with $\|q^{(0)}\|_{L^\infty(Q;\bbR^d)}\le R$. The policy iteration method proceeds by iterating on $k\ge 0$ the following steps \cite{CaCaGo-COCV21}:
\begin{enumerate}[label=(\roman*)]
	\item Solve the following linear FP equation for $m^{(k)}$, given the current policy $q^{(k)}$:%
	\begin{equation}
	\left\{\begin{array}{ll}
	\partial_t m^{(k)} - \eps\Delta m^{(k)} - {\rm div}(m^{(k)} q^{(k)}) = 0\, &\text{in } Q,\\
	m^{(k)}(\bx,0)=m_0(\bx)\, &\text{in } \bbT^d.
	\end{array}\right.
	\end{equation}
	
	\item Solve the following linear equation for $u^{(k)}$, given $q^{(k)}$ and $m^{(k)}$:
	\begin{equation}
	\left\{\begin{array}{@{}l@{\;}l@{}}
	-\partial_t u^{(k)}-\eps\Delta u^{(k)} + q^{(k)}\cdot \Du^{(k)} - L(q^{(k)}) = b(\bx) +  F( m^{(k)})\, &\text{in } Q,\\ 
	u^{(k)}(\bx,T) = u_T(\bx)\, &\text{in } \bbT^d. 
	\end{array}\right.
	\end{equation} 
	
	\item Update the policy
	\begin{equation}
	q^{(k+1)}(\bx,t) = \argmax_{|q|\le R} \left[q\cdot \Du^{(k)} - L(q)\right].
	\end{equation}
	
\end{enumerate}
The policy function $q^{(k)}: Q\to \bbR^d$ represents the optimal strategy of agents at the $k$th iteration. It decouples the system \eqref{EQ:MFG} into several iterations of linear PDE solves (step (i) and (ii)), and optimization problems (step (iii), can be solved point-wisely) with closed-form solutions in some cases. 
Because of the decoupling of the strongly coupled forward-backward systems, the PDEs in steps (i) and (ii) can be solved using classical numerical methods that integrate in one direction. Additionally, introducing the policy function linearizes both equations, converting the nonlinearity of MFGs into a fixed-point iteration of linear PDEs. This ensures the efficiency of policy iteration. Together with theoretical analysis guaranteeing its linear convergence \cite{CiGiMa-DGA20,CaTa-JMAA22},
this provides us with a solid foundation for developing methods for inverse MFGs.

\section{Policy iteration for inverse MFG problem\label{sec:alg}}

In this section, we propose a policy iteration method for inverse MFGs introduced at the end of \cref{sec:MFG}, inspired by the policy iteration method for solving MFG equations discussed in \cref{sec:Policy-MFG}. Before detailing the proposed method, we highlight the drawbacks of directly applying policy iteration to a least-squares method in \cref{sec:direct-LS}. We then introduce our tailored policy iteration method for inverse MFGs in \cref{sec:Policy-IP-MFG}.

\subsection{Direct application to least-squares of data misfit\label{sec:direct-LS}}
The standard approach to solving inverse problems is to use the \textbf{direct least-squares (LS) method}, formulating it as a PDE-constrained optimization problem to minimize the squared $L^2$-misfit between the observed data $\DataFun(\bx)$ and the corresponding measurement $\MeasOpr u$ of the solutions of the MFG model \eqref{EQ:MFG}:
\begin{equation}
\label{eq:direct-LS}
\begin{aligned}
\underset{b(\bx)}{\text{minimize}}\;\; & \Phi(b) := \frac{1}{2} \int_{\bbT^d} [\MeasOpr u(\bx) - \DataFun(\bx)]^2\, d\bx \\
 \text{subject to}\;\; & (u,m) \text{ solves the MFG model \eqref{EQ:MFG}}.
\end{aligned}
\end{equation}
To solve this PDE-constrained optimization, gradient information is required, which can be evaluated using adjoint methods \cite{ReCa-book15,HiPiUl-book09,BoSc-book11,Tr-AMS10}.
Details of the application for the two cases: (i) $ \MeasOpr u:=u(\bx,0)$ and (ii) $\MeasOpr u:=\partial_t u(\bx, T)$ 
can be found in \cref{sec:LS-U0} and \cref{sec:LS-UtN} 
respectively. Because of the coupled structure of the state equation \eqref{EQ:MFG}, the adjoint equations (such as \eqref{eq:adj-LS-U0} and \eqref{eq:adj-LS-UtN}) also form strongly coupled systems of forward and backward equations, which need to be solved iteratively (similar to policy iteration for solving MFGs \eqref{EQ:MFG}; note that policy is not introduced for adjoint equations because of its linearity). However, this direct LS method faces challenges both in terms of optimization and computation. The optimization problem \eqref{eq:direct-LS} is nonconvex and PDE-constrained, thus its gradient may only guide toward a local optimizer, and its performance highly depends on the initialization. In addition, each gradient evaluation requires solving two forward-backward coupled PDE systems, meaning that computing the full policy iterations for solving the MFG and similar iterative methods for the adjoint equation is required for every gradient evaluation, making the process computationally intensive.

\subsection{Proposed algorithm: policy iteration method for inverse MFGs\label{sec:Policy-IP-MFG}}

Because of these issues, we aim to design a more effective algorithm for inverse MFG problems than the direct LS method. Specifically, we seek to develop a policy iteration method for inverse MFG problems. Note that the observed data depends only on the value function $u$, which appears only in step (ii) of the policy iteration method. This motivates us to incorporate the inversion step exclusively in step (ii) and formulate the following \textbf{policy iteration method for inverse MFG problems}, which requires only a few iterations of solving linear PDEs and linear inverse problems.

Choose an initial vector field $q^{(0)}\in C^{1,0}(Q;\bbR^d)$.
The proposed \textbf{policy iteration method for inverse MFGs} proceeds by iterating on $k\ge 0$ the following three steps:
\begin{enumerate}[label=(\roman*)]
    \item Solve the following linear FP equation for $m^{(k)}$, given the current policy $q^{(k)}$:
    \begin{equation}\label{eq:FP-policy}
    \left\{\begin{array}{ll}
        \partial_t m^{(k)} - \eps\Delta m^{(k)} - {\rm div}(m^{(k)} q^{(k)}) = 0 \,&\text{in } Q,\\
        m^{(k)}(\bx,0)=m_0(\bx)\, &\text{in }  \bbT^d.
    \end{array}\right.
    \end{equation}

    \item\label{IT:Pol Eval Inv} Solve the following linear inverse problem: determine $b^{(k)}(\bx)$ such that $\MeasOpr u^{(k)}\approx \DataFun(\bx)$ with  $u^{(k)}$ the solution to the following linear PDE, given $q^{(k)}$ and $m^{(k)}$:
    \begin{equation}\label{eq:linear-HJB}
    \left\{\begin{array}{@{}l@{\,}l@{}}
        -\partial_t u^{(k)}-\eps\Delta u^{(k)} + q^{(k)}\cdot \Du^{(k)} - L(q^{(k)}) = b^{(k)} + F(m^{(k)})\, &\text{in } Q,\\
        u^{(k)}(\bx,T) = u_T(\bx)\,&\text{in }  \bbT^d.
    \end{array}\right.
    \end{equation} 

    \item  Update the policy: 
    \begin{equation}
        q^{(k+1)}(\bx,t) = \argmax_q %
        \left[q\cdot \Du^{(k)}- L(q)\right].
    \end{equation}

\end{enumerate}

The step (iii) can be further rewritten
using the properties of Legendre transform, 
\begin{equation}\label{eq:q-DpHDu}
q^{(k+1)}(\bx,t) = \argmax_q %
\left[q\cdot \Du^{(k)}- L(q)\right] = \DpH(\Du^{(k)}).
\end{equation}
For example, for a quadratic Hamiltonian $H(p)= \frac{1}{2}|p|^2$, the corresponding Lagrangian is $L(q)=\frac{1}{2}|q|^2$, and the updated policy is $q^{(k+1)}(\bx,t) =\Du^{(k)}(\bx,t)$.

Compared with the policy iteration method for MFGs, steps (i) and (iii) remain unchanged, while step (ii) is changed from solving a linear PDE to solving the inverse problem with this linear PDE as a constraint.
Furthermore, 
the solution $u^{(k)}$ in step (ii) depends linearly on the obstacle function $b^{(k)}(\bx)$, meaning that the observation measurement $\MeasOpr u^{(k)}$ is also a linear measurement of $b^{(k+1)}$. Therefore, the inverse problem of reconstructing $b^{(k)}$ such that $\MeasOpr u^{(k)}\approx \DataFun(\bx)$ is a linear inverse problem.

This \textbf{linear parabolic inverse source problem} in step (ii) can be solved using the standard linear least-squares method. Similar to the direct LS method, we can formulate the step (ii) as a PDE-constrained optimization problem
\begin{equation}
\label{eq:LS-HJB}
\begin{aligned}
\underset{b^{(k)}(\bx)}{\text{minimize}}\;\; & \Phi(b^{(k)}) = \frac{1}{2} \int_{\bbT^d} [\MeasOpr u^{(k)}(\bx) - \DataFun(\bx)]^2\, d\bx \\
\text{subject to}\;\; & u^{(k)} \text{ solves the linear PDE \eqref{eq:linear-HJB}}.
\end{aligned}
\end{equation}
This is a convex optimization problem with a quadratic objective and a linear constraint, meaning that its local minimizers are also global minimizers. It can be solved using gradient-based algorithms (such as quasi-Newton methods, for example, BFGS), with gradients evaluated through adjoint methods (details in \cref{sec:Adj-Policy-U0}). 
Since its optimality condition can be viewed as a linear system, conjugate gradient, and other Krylov methods can be applied to further accelerate the optimization.
Because of the convexity of the problem, the initialization of the optimization algorithm is not a critical factor and can be chosen as the optimizer $b^{(k-1)}(\bx)$ from the previous policy iteration. 
In addition to convexity, another advantage is the reduced computational cost. Each gradient evaluation requires solving only two linear PDEs, compared to the direct LS method, which requires solving two coupled PDE systems. 

The efficiency of solving the linear inverse problem in step (ii) can be further enhanced in certain cases. For example,  when the data $\DataFun(\bx)$ are observations of $\partial_t u$ at the final time (case ii), i.e., $\MeasOpr u := \partial_t u(\bx, T)$, the optimization problem \eqref{eq:LS-HJB} 
can be solved in a one-shot approach by evaluating \eqref{eq:linear-HJB} at the final time:
 $b^{(k)}(\bx)  = \cB(q^{(k)},m^{(k)})(\bx) $,
where $\cB$ is the solution operator of the inverse problem in step (ii), and has a closed-form formula in this case,
\begin{equation}\label{eq:B-LS-HJB-UtN}
\cB(q,m)(\bx) = -\DataFun(\bx) - \eps\Delta u_T + q(\bx,T)\cdot \Du_T -L(q(\bx,T)) - F(m(\bx,T)).
\end{equation}
This discussion demonstrates the solvability of the inverse problem in step (ii), and reformulates it back into solving the following linear PDE for $u^{(k)}$, given $m^{(k)},q^{(k)}$,
\begin{equation}\label{eq:linear-HJB-B}
\left\{\begin{array}{@{}l@{\;}l@{}}
-\partial_t u^{(k)}-\eps\Delta u^{(k)} + q^{(k)}\cdot \Du^{(k)} -  L(q^{(k)}) = \cB(q^{(k)},m^{(k)})(\bx)+ F(m^{(k)})\, &\text{in } Q,\\
u^{(k)}(\bx,T) = u_T(\bx)\,&\text{in }  \bbT^d.
\end{array}\right.
\end{equation} 
The solution $u^{(k)}$ satisfies $\partial_t u^{(k)}(\bx, T) = \DataFun(\bx)$. 
Moreover, 
using the fact that $q^{(k)} (\bx, T)= \DpH(\Du^{(k-1)}(\bx, T) )=  \DpH(\Du_T (\bx) )$ is the optimizer for Legendre transform $H (\Du_T) = \max_q [q\cdot \Du_T - L(q)]$,
$\cB(q^{(k)},m^{(k)})(\bx)$ in this case can be simplified as 
\begin{equation}\label{eq:B-simplify}
\cB(q^{(k)},m^{(k)})(\bx) = -\DataFun(\bx) - \eps\Delta u_T + H (\Du_T) - F(m^{(k)}(\bx,T)).
\end{equation}

\section{Convergence of the policy iteration method for inverse MFG problems}

The policy iteration method for MFGs was proven to have a linear convergence rate in ~\cite{CiGiMa-DGA20,CaTa-JMAA22}.
Although our proposed policy iteration method for inverse MFG problems has a similar structure, the aforesaid results do not directly cover our method because of the additional inversion in step (ii).
In this section, we present a convergence result for the case when the inverse policy iteration method is applied to $\partial_t u(\bx,T)$ data and the time horizon is sufficiently small. 
The policy iteration method for inverse MFGs converges uniformly, which we prove by showing it is a contractive fixed-point iteration in \cref{sec:converge-uniform}.
Furthermore, we prove its rate of convergence is linear in \cref{sec:converge-linear}. All the proofs in this section are provided for a quadratic Hamiltonian,  %
but can be extended to more general cases.

\subsection{Notations}

Before discussing the theoretical framework, we first introduce the notations for spaces and norms (adapted from \cite{LaSoTa-AMO23}) used in the theorems and proofs. 

The vector norm $|q|$ and the inner product $p\cdot q$ are in the standard form, for $p, q\in \bbR^d$.
For $1\leq r \leq \infty$, $L^r(Q)$ is the usual Lebesgue space with norm $\|\cdot\|_{L^r(Q)}$, and $L^r(Q; \bbR^d)$ and $L^\infty(Q; \bbR^d)$ are the corresponding vector-valued Lebesgue space. $W_r^k$ is the standard Sobolev spaces with weak derivatives up to the order of $k$ in $L^r$-norm. The anisotropic Sobolev space $W_r^{2,1}(Q)$ includes function $u$ such that $\partial_t^\delta \nabla^\sigma u\in L^r(Q)$ for all $|\sigma|+ 2\delta \leq 2$, endowed with norm 
$
\|u\|_{W^{2,1}_r(Q)}=(\int_{Q}\sum_{|\sigma|+2{\delta}\leq 2} |\partial_t^{{\delta}}\nabla^{\sigma} u |^r \, d\bx dt)^{\frac1r}.
$ The trace of $W^{2,1}_r(Q)$ is given by the fractional Sobolev class $W^{2-\frac{2}{r}}_r({\bbT^d})$.
$C^{1,0}(Q)$ is the space of continuous functions on $Q$ with continuous derivatives in the $x$-variable, endowed with the norm $|u|^{(1)}_Q:=\| u \|_{L^\infty(Q)} + \|\Du\|_{L^\infty(Q)} $.
For $0< \alpha < 1$, the H\"older space $C^{\alpha,{\alpha}/{2}}(Q)$ is endowed with the norm
\begin{equation}\label{eq:Holder-norm-alpha}
| u|^{(\alpha)}_Q=\| u \|_{L^\infty(Q)}+\sup_{(\bx,t)\neq(\by,s)\in Q}\frac{| u(\bx,t)-u(\by,s)| }{(d(\bx,\by)^2+| t-s| )^{\frac{\alpha}{2}}},
\end{equation}
where $d(\bx,\by)$ stands for the geodesic distance from $\bx$ to $\by$ in ${\bbT^d}$. The H\"older space $C^{1+\alpha,{(1+\alpha)}/{2}}(Q)$ is endowed with the norm
\begin{equation}\label{eq:Holder-norm-alpha+1}
| u|^{(1+\alpha)}_Q=\| u \|_{L^\infty(Q)} +\sum_{i=1}^d| \partial_{x_i}u|^{(\alpha)}_Q+\sup_{(\bx,t)\neq(\by,s)\in Q}\frac{| u(\bx,t)-u(\by,s)| }{| t-s| ^{\frac{1+\alpha}{2}}} .
\end{equation}

\subsection{Uniform convergence theorem\label{sec:converge-uniform}}

First, we provide a uniform convergence theorem for the proposed policy iteration method for inverse MFG problems.

\begin{theorem}\label{THM:Convergence} 
	Under the assumptions:
	\begin{enumerate}[leftmargin=*, label=($\cA$-\roman*)]
		\item \label{ASS:Init Term Cond} The initial and terminal conditions $u_T \in W^2_{\infty}(\bbT^d)$, $m_0\in W^2_r(\bbT^d)$ for some $r>d+2$, $m_0\ge \underline{m}>0$ for some constant $\underline{m}$, and $\int_{\bbT^d} m_0(\bx)\,d\bx = 1$.
		\item \label{ASS:F} The interaction cost function $F:\bbR_{\ge 0}\rightarrow \bbR$ is strictly increasing and locally Lipschitz continuous. 
		\item \label{ASS:Data} Data $\DataFun\in L^{\infty}(\bbT^d)$. 
		\item \label{ASS:q0} Initial policy $q^{(0)}\in C^{1,0}(Q;\bbR^d)$.
	\end{enumerate}
	There exists a $\bar{T}$ such that  $\forall T\in (0,\bar{T}]$, 
 the sequence $\{b^{(k)}\}_{k\ge 0}$, generated by the policy iteration method for the inverse MFG with $\partial_t u(\bx, T)$ data and a quadratic Hamiltonian $H(p)=\frac{1}{2}|p|^2$, converges uniformly on $\bbT^d$ to a solution $b^*$ of the inverse problem, i.e.,%
\begin{equation}\label{eq:converge-b-fix-point}
\lim_{k\to\infty} \|b^{(k)}(\bx) - b^*(\bx)\|_{L^\infty(\bbT^d)} = 0.
\end{equation}
\end{theorem}

\noindent\emph{Remark.} The value of $\bar{T}$ depends on the given information of the inverse MFG problems: the diffusion coefficient $\eps$, the terminal value function $u_T$, the initial density function $m_0$, the interaction cost function $F$, and the measured data $\DataFun$. It also depends on the initialization $q^{(0)}$ for the policy iteration method.

\edit{\noindent\emph{Remark.} The monotonicity of $F$ required in the assumption \ref{ASS:F} is not used in the following convergence proofs for the proposed algorithm in solving the inverse MFG problem. However, the monotonicity of $F$ is required to ensure the well-posedness of the forward MFG problem \cite{LaLi-JJM07}.}

\begin{proof}[Proof of \Cref{THM:Convergence}] 
	
	The proof closely follows \cite{LaSoTa-AMO23}  utilizing the contraction mapping theorem. 
Specifically, 
	we use several ``constants" for the embedding estimates in the proof (typically using the notation $C$ with subscription/superscription). These ``constants"
depend only on the given information of the inverse MFG problems: $\eps$,  $u_T$,  $m_0$, $F$, and $\DataFun$. Occasionally, these ``constants" also depend on the final time $T$, but remain bounded for bounded values of $T$.

From the discussion in \cref{sec:Policy-IP-MFG}, the policy iteration method for inverse MFGs reduces to iterations of solving linear PDEs \eqref{eq:FP-policy} and \eqref{eq:linear-HJB-B}. In the quadratic Hamiltonian setting with the policy replaced by $q^{(k)} = \Du^{(k-1)}$, it 
can be viewed as a map $\cT$ from $(u^{(k-1)}, m^{(k-1)})$ to $(u^{(k)}, m^{(k)})$, i.e., $(u^{(k)}, m^{(k)}) = \cT (u^{(k-1)}, m^{(k-1)})$, 
with the operator defined as $\cT: (u,m) \mapsto (\hat{u}, \hat{m})$ with $(\hat{u}, \hat{m})$ the solution to 
\begin{equation}\label{EQ:T}
\left\{\begin{array}{ll}
\partial_t \hat{m} - \eps\Delta \hat{m} - {\rm div}(\hat{m} \Du) = 0\, &\text{in } Q,\\
-\partial_t \hat{u}-\eps\Delta \hat{u} + \Du\cdot \Dhu - \frac12|\Du|^2 = \cB(\Du,\hat{m})(\bx) +F(\hat{m})\, &\text{in } Q,\\
\hat{u}(\bx,T) = u_T(\bx),\ \hat{m}(\bx,0)=m_0(\bx)\, &\text{in } \bbT^d,
\end{array}\right.
\end{equation}
where $\cB$ is defined in \eqref{eq:B-LS-HJB-UtN} with $L(q)=\frac{1}{2}|q|^2$ (or, for simplicity, we can instead define $\cB$ as in \eqref{eq:B-simplify}, the following proofs still hold).
The solution to the linear inverse problem in step (ii)  is 
$b^{(k)}(\bx) = \cB(q^{(k)}, m^{(k)})= \cB(\Du^{(k-1)}, m^{(k)})$.

The main idea is to show $\cT$ is a contraction map, thus the sequence $\{ (u^{(k)}, m^{(k)}) \}_{k\ge 0}$ converges to a fixed point $(u^*, m^*)$ corresponding to a solution 
$b^*:=\cB(q^{*},m^{*})$ 
to the inverse problem. To achieve this, we define the following space 
as
the domain of $\cT$:
\begin{equation*}
    \XTM := \{(u,m): u\in C^{1,0}(Q)\cap W^{2,1}_r(Q), m\in C^{1,0}(Q), \|u\|_{W^{2,1}_r(Q)} + |u|^{(1)}_Q + |m|^{(1)}_Q \le M_1\},
\end{equation*}
with the constant upper bound $M_1 :=  2(|m_0|^{(1)}_{\bbT^d} + |u_T|^{(1)}_{\bbT^d})$.

In order to prove the contraction property of operator $\cT$, we need the boundness of the coefficients and source terms of \eqref{EQ:T} for applying parabolic estimates \Cref{PROP:Parabolic Estimate Strong,PROP:Parabolic Estimate Weak}.
To achieve this, we utilize a related operator $\cTK$,
defined as $\cTK: (u,m)\mapsto(\hat{u},\hat{m})$  with $(\hat{u},\hat{m})$ the solution to 
\begin{equation}\label{EQ:T_K}
\left\{\begin{array}{ll}
\partial_t \hat{m} - \eps\Delta \hat{m} - \psi(\Du)\cdot \Dhm  - (\Delta u)\hat{m} = 0\, &\text{in } Q,\\
-\partial_t \hat{u}-\eps\Delta \hat{u} + \psi(\Du)\cdot \Dhu - \frac12|\psi(\Du)|^2 =
 \cB(\psi(\Du),\varphi(\hat{m}))(\bx) 
+F(\varphi(\hat{m}))\, &\text{in } Q,\\
\hat{u}(\bx,T) = u_T(\bx),\ \hat{m}(\bx,0)=m_0(\bx)\, &\text{in } \bbT^d,
\end{array}\right.
\end{equation} 
where $\Du$ and $\hat{m}$ are regularized by two bounded, globally Lipschitz functions: $\varphi:\bbR\rightarrow \bbR$ such that $\varphi(z) = z$ for all $z\in [1/K, K]$ and $\varphi(z)\in [1/(2K), 2K]$ for all $z\in \bbR$, and $\psi:\bbR^d\rightarrow \bbR^d$ such that $\psi(p) = p$ for all $|p|\le K$ and $|\psi(p)| \le 2K$ for all $p\in \bbR^d$, with the constant $K$ defined as
\begin{equation}\label{eq:K}
K :=  \max\{ {2}{\underline{m}}^{-1},\, 2|m_0|^{(1)}_{\bbT^d}, \, 2|u_T|^{(1)}_{\bbT^d} \}.
\end{equation}
Since $F$ is locally Lipschitz continuous and $\psi$ and $\varphi$ are bounded, globally Lipschitz continuous, their compositions are also bounded and globally Lipschitz continuous, i.e., there exist some  constants $C_0$ and $C_0'$
depending only on the problem information and $K$, for all $(u_1,m_1), (u_2,m_2)\in \XTM$, 
\begin{equation}\label{eq:Lipschitz-Tk}
\begin{array}{r@{\;}l}
\left\|\left|\frac12|\psi(\Du_1)|^2 - \frac12|\psi(\Du_2)|^2 \right| + |\psi(\Du_1) - \psi(\Du_2)|
\right\|_{L^\infty(Q)} & \le  C_0 |u_1-u_2|^{(1)}_Q, \\
\left\|F(\varphi(m_1)) - F(\varphi(m_2))\right\|_{L^\infty(Q)} & \le C_0 |m_1-m_2|^{(1)}_Q,\\
\left\|\cB(\psi(\Du_1),\varphi(m_1))(\bx) - \cB(\psi(\Du_2),\varphi(m_2))(\bx)\right\|_{L^\infty(\bbT^d)}  
& \le C_0(|u_1-u_2|^{(1)}_Q + |m_1-m_2|^{(1)}_Q), 
\end{array}
\end{equation}
and for all $(u,m)\in \XTM$,
using the assumptions $u_T\in W^2_{\infty}(\bbT^d)$, and $\DataFun\in L^{\infty}(\bbT^d)$,
\begin{equation}\label{eq:bound-Tk}
\begin{array}{rl}
\left\|\frac12|\psi(\Du)|^2 + |\psi(\Du)| + |F(\varphi(m))|\right\|_{L^\infty(Q)}  & \le   C_0',\\
\left\|\cB(\psi(\Du),\varphi(\hat{m}))(\bx) \right\|_{L^\infty(\bbT^d)}  & \le    C_0'.
\end{array}
\end{equation}

The following proofs can be separated into four steps: Steps 1 and 2 prove the map $\cTK$ maps $\XTM$ to itself and the map is contractive; Step 3 shows $\cT=\cTK$ in a subset of $\XTM$ and maintains the same properties; Step 4 uses the contraction mapping theorem to prove the sequence  
$b^{(k)} = \cB(q^{(k)},m^{(k)})$ 
converges to the solution of the inverse MFG problem.

\emph{Step 1: $\cTK$ maps $\XTM$ into itself.} 
$\forall (u,m)\in \XTM$, we want to show that $ (\hat{u},\hat{m}):=\cTK(u,m)  \in \XTM$, where $\cTK$ is defined in \eqref{EQ:T_K}. 
The proof below applies embedding estimates to $\hat{m}$ and $\hat{u}$ respectively as the solution to 
\eqref{EQ:T_K}, based on \Cref{PROP:Parabolic Estimate Weak,PROP:Parabolic Estimate Strong} for linear parabolic equations.

First, we consider the solution $\hat{m}$ to the first equation in \eqref{EQ:T_K} with coefficients $\Delta u$ and $\psi(\Du)$. 
From the definition $u\in\XTM$, $ \|\Delta u\|_{L^r(Q)} \le M_1$. From \eqref{eq:bound-Tk} and H\"{o}lder's inequality,
\begin{equation}\label{eq:bound-psi-Du}
   \|\psi(\Du)\|_{L^{\infty}(Q; \bbR^d)}\le C_0', \qquad \|\psi(\Du)\|_{L^r(Q; \bbR^d)} \le T^{1/r} C_0'.
\end{equation}
Using the parabolic estimate \Cref{PROP:Parabolic Estimate Weak}, we have 
$\|\hat{m}\|_{W^{2,1}_r(Q)} \le C_1.$
The embedding result of \Cref{PROP:Embed Holder Sobolev} then implies 
\begin{equation}\label{eq:bound-hatm-Qnorm-long}
    |\hat{m}|^{(2-\frac{d+2}{r})}_Q \le C_2(\|\hat{m}\|_{W^{2,1}_r(Q)} + \|m_0\|_{W^{2-\frac2r}_r(\bbT^d)})\le C_2',
\end{equation}
combining with \Cref{LEM:Embed 1 1+alpha} provides the estimate
\begin{equation}\label{eq:bound-hatm-Qnorm}
    |\hat{m}|^{(1)}_Q \le |m_0|^{(1)}_{\bbT^d} + T^{\frac12 - \frac{d+2}{2r}} C_2'.
\end{equation}

Secondly, consider the solution $\hat{u}$ in the second equation of \eqref{EQ:T_K}, the coefficient $\psi(\Du)\in L^{\infty}(Q;\bbR^d)$ from \eqref{eq:bound-psi-Du}.  From \eqref{eq:bound-Tk}, the source term $\frac12|\psi(\Du)|^2 + \cB(\psi(\Du), \varphi(\hat{m})) + F(\varphi(\hat{m}))$ is bounded by $2C_0'$
in $L^{\infty}(Q)$, thus it is also bounded by $2T^{\frac{1}{2r}}C_0'$
in $L^{2r}(Q)$ 
(using H\"{o}lder's inequality).
Therefore, \Cref{PROP:Parabolic Estimate Strong} provides the following estimate 
\begin{equation*}
    \|\hat{u}\|_{W^{2,1}_{2r}(Q)} \le C_3(2T^{\frac{1}{2r}}C_0' + \|u_T\|_{W^{2-\frac1r}_{2r}(\bbT^d)}) \le C_3',
\end{equation*}
and a further application of \Cref{LEM:Embed r 2r} yields the estimate of $\hat{u}$ in $W^{2,1}_{r}(Q)$,
\begin{equation}\label{eq:bound-hatu-Wnorm}
    \|\hat{u}\|_{W^{2,1}_{r}(Q)} \le T^{\frac{1}{2r}} \|\hat{u}\|_{W^{2,1}_{2r}(Q)} \le T^{\frac{1}{2r}}C_3'.
\end{equation}
Following the same derivation in \eqref{eq:bound-hatm-Qnorm-long}, \Cref{PROP:Embed Holder Sobolev} implies
$
    |\hat{u}|^{(2-\frac{d+2}{r})}_Q 
    \le C_3''.
    $
This together with \Cref{LEM:Embed 1 1+alpha} provides the estimate of $\hat{u}$ in $C^{1,0}(Q)$,
\begin{equation}\label{eq:bound-hatu-Qnorm}
    |\hat{u}|^{(1)}_Q \le |u_T|^{(1)}_{\bbT^d} + T^{\frac12 - \frac{d+2}{2r}} C_3''.
\end{equation}

Combining the estimates \eqref{eq:bound-hatm-Qnorm}, \eqref{eq:bound-hatu-Wnorm} and \eqref{eq:bound-hatu-Qnorm}, we obtain
\begin{align*}
    \|\hat{u}\|_{W^{2,1}_r(Q)} + |\hat{u}|^{(1)}_Q + |\hat{m}|^{(1)}_Q &\le T^{\frac{1}{2r}}C_3' + |u_T|^{(1)}_{\bbT^d} + T^{\frac12 - \frac{d+2}{2r}}C_3'' + |m_0|^{(1)}_{\bbT^d} + T^{\frac12 - \frac{d+2}{2r}} C_2'\\
    &\le T^{\frac{1}{2r}}C_3' + T^{\frac12 - \frac{d+2}{2r}} ( C_2' + C_3'') + \tfrac{1}{2}M_1,
\end{align*}
which can be enforced to be less than $M_1$ for a sufficiently small $T$, since the exponent on $T$ is positive.
This shows that $ \cTK(u,m) = (\hat{u},\hat{m})  \in \XTM$, as desired.

\emph{Step 2: $\cTK: \XTM\rightarrow \XTM$ is a contraction operator.} $\forall (u_1,m_1), (u_2,m_2)\in \XTM$, we want to show that there exists a ``constant" $0<\Gamma<1$ (independent of the choices of $(u_1,m_1)$ and $(u_2,m_2)$, but can depend on the inverse problem information or $T$, and stay bounded for bounded values of $T$), such that $\cTK(u_1,m_1) = (\hat{u}_1, \hat{m}_1)$ and $\cTK(u_2,m_2) = (\hat{u}_2, \hat{m}_2)$ satisfy 
\begin{equation}\label{eq:def-contraction}
    \|\hat{u}_1 - \hat{u}_2\|_{W^{2,1}_r(Q)} + |\hat{u}_1 - \hat{u}_2|^{(1)}_Q + |\hat{m}_1 - \hat{m}_2|^{(1)}_Q
\le \Gamma(\|u_1 - u_2\|_{W^{2,1}_r(Q)} + |u_1 - u_2|^{(1)}_Q + |m_1 - m_2|^{(1)}_Q).
\end{equation}
The following proofs are also based on parabolic estimates \Cref{PROP:Parabolic Estimate Weak,PROP:Parabolic Estimate Strong}, 
by applying to the difference of two solutions $\barU := \hat{u}_1 - \hat{u}_2$ and $\barM := \hat{m}_1 - \hat{m}_2$. 

Firstly, we obtain a linear parabolic equation for $\barM$ by a subtraction of the first equation in \eqref{EQ:T_K} for each pair of $(\hat{m}_1, u_1)$ and $(\hat{m}_2, u_2)$,
\begin{equation}\label{eq:barM-PDE}
\left\lbrace \begin{array}{@{}l}
    \partial_t \barM - \eps\Delta \barM - \psi(\Du_1)\cdot \DbM + (\Delta u_1)\barM
= [\psi(\Du_1) - \psi(\Du_2)]\cdot \Dhm _2 - \Delta(u_1 - u_2)\hat{m}_2,\\
\barM(\bx,0) = 0.
\end{array}\right. 
\end{equation}
Its coefficients $\psi(\Du_1)$ and $\Delta u_1$ are bounded in $L^{r}$ from the same derivations \eqref{eq:bound-psi-Du} in Step 1. 
Its source term is also bounded in $L^r(Q)$, by H\"{o}lder's inequality, together with the definition $\hat{m}_2\in \XTM$ and Lipschitz regularity \eqref{eq:Lipschitz-Tk}, i.e.,
\begin{align*}
  &  \|[\psi(\Du_1) - \psi(\Du_2)]\cdot \Dhm _2\|_{L^r(Q)} 
    \le T^{\frac1r} \|[\psi(\Du_1) - \psi(\Du_2)]\cdot \Dhm _2\|_{L^\infty(Q)} \\
   &\qquad \leq T^{\frac1r} \||\psi(\Du_1) - \psi(\Du_2)|\|_{L^\infty(Q)}\| |\Dhm _2| \|_{L^\infty(Q)}
    \le T^{\frac1r} C_0 M_1 |u_1 - u_2|^{(1)}_Q,\\
    & \|\Delta(u_1 - u_2)\hat{m}_2\|_{L^r(Q)} \le \|\Delta(u_1 - u_2)\|_{L^r(Q)}\|\hat{m}_2\|_{L^{\infty}(Q)} \le M_1 \|u_1 - u_2\|_{W^{2,1}_r(Q)}.
\end{align*}
Therefore, \Cref{PROP:Parabolic Estimate Weak} on \eqref{eq:barM-PDE} informs that for sufficiently small $T$, the following holds
\begin{align*}
    \|\barM\|_{W^{2,1}_r(Q)} \le C_4(|u_1 - u_2|^{(1)}_Q + \|u_1 - u_2\|_{W^{2,1}_r(Q)}),
\end{align*}
with a further application of \Cref{LEM:Embed Holder Sobolev Cor} providing the estimate for $\barM$ in $C^{1,0}(Q)$,
\begin{equation}\label{eq:bound-barM-Qnorm}
    |\barM|^{(1)}_Q \le T^{\frac12 - \frac{d+2}{2r}}C_4'(|u_1 - u_2|^{(1)}_Q + \|u_1 - u_2\|_{W^{2,1}_r(Q)}). 
\end{equation}

Secondly, we obtain the linear parabolic equation for $\barU$ by a subtraction of the second equation of \eqref{EQ:T_K} for each triple of $(\hat{u}_1, u_1, \hat{m}_1)$ and $(\hat{u}_2, u_2, \hat{m}_2)$ provides
\begin{equation}\label{eq:barU-PDE}
\left\lbrace \begin{array}{@{}l@{\;}l}
-\partial_t\barU - \eps\Delta\barU + \psi(\Du_1)\cdot \DbU =& -\left[\psi(\Du_1) - \psi(\Du_2)\right]\cdot \Dhu_2 + \left[\frac12|\psi(\Du_1)|^2 - \frac12|\psi(\Du_2)|^2\right]\\
&+\cB(\psi(\Du_1),\varphi(\hat{m}_1))(\bx) - \cB(\psi(\Du_2),\varphi(\hat{m}_2)) (\bx)\\
&+ F(\varphi(\hat{m}_1)) - F(\varphi(\hat{m}_2)),\\
\barU(\bx,T) = 0.&
\end{array}\right. 
\end{equation}
Its coefficient $\psi(\Du_1)$ is bounded in $L^{\infty}(Q; \bbR^d)$ from \eqref{eq:bound-psi-Du}.
Similarly as the derivation for $\barM$, 
the $L^{\infty}(Q)$ norm of the source term in \eqref{eq:barU-PDE} is bounded by $C_5(|u_1-u_2|^{(1)}_Q + |\hat{m}_1-\hat{m}_2|^{(1)}_Q)$ from the Lipschitz properties \eqref{eq:Lipschitz-Tk}, thus
its $L^r(Q)$ norm can be bounded by 
$T^{\frac1r}C_5(|u_1-u_2|^{(1)}_Q + |\hat{m}_1-\hat{m}_2|^{(1)}_Q)$ from H\"{o}lder's inequality.
Therefore, an application of \Cref{PROP:Parabolic Estimate Strong} to \eqref{eq:barU-PDE} provides the estimate of $\barU$ in $W^{2,1}_{r}(Q)$
\begin{equation}\label{eq:bound-barU-Wnorm}
\|\barU\|_{W^{2,1}_r(Q)} \le T^{\frac1r}C_5(|u_1-u_2|^{(1)}_Q + |\hat{m}_1-\hat{m}_2|^{(1)}_Q).
\end{equation}
Moreover, from \Cref{LEM:Embed Holder Sobolev Cor}, we obtain the estimate of $\barU$ in $C^{1,0}(Q)$,
\begin{equation}\label{eq:bound-barU-Qnorm}
|\barU|^{(1)}_Q \le T^{\frac12 - \frac{d+2}{2r}}C_5'(|u_1-u_2|^{(1)}_Q + |\hat{m}_1-\hat{m}_2|^{(1)}_Q).
\end{equation}
Finally, we combine the estimates \eqref{eq:bound-barM-Qnorm}, \eqref{eq:bound-barU-Wnorm} and \eqref{eq:bound-barU-Qnorm} and conclude that 
\begin{align*}
    &\|\hat{u}_1 - \hat{u}_2\|_{W^{2,1}_r(Q)} + |\hat{u}_1 - \hat{u}_2|^{(1)}_Q + |\hat{m}_1 - \hat{m}_2|^{(1)}_Q = \|\barU\|_{W^{2,1}_r(Q)} + |\barU|^{(1)}_Q + |\barM|^{(1)}_Q\\
    & \leq (T^{\frac1r}C_5 + T^{\frac12 - \frac{d+2}{2r}}(C_4'+C_5'))(\|u_1 - u_2\|_{W^{2,1}_r(Q)} + |u_1 - u_2|^{(1)}_Q + |m_1 - m_2|^{(1)}_Q).
\end{align*}
Define $\Gamma:=T^{\frac1r}C_5 + T^{\frac12 - \frac{d+2}{2r}}(C_4'+C_5')$. Since the exponent on $T$ is positive, the factor $\Gamma$ can be enforced to be less than $1$ when $T$ is sufficiently small. This meets our goal of \eqref{eq:def-contraction}, meaning that $\cTK: \XTM\rightarrow \XTM$ is a contraction map.

\emph{Step 3: $\cT$ is a contraction map on the space $\XTM\cap \{1/K\le m\le K,\ |u|^{(1)}_Q \le K\}$.} 
We have just shown that $\cTK$ is a contraction map on $\XTM$ (Step 2). If we can additionally show that $\cTK$ maps the space $\XTM\cap \{1/K\le m\le K,\ |u|^{(1)}_Q \le K\}$ to itself, then comparing \eqref{EQ:T} and \eqref{EQ:T_K} shows that $\cT=\cTK$ on this space (and is a contraction). Therefore, we want to show
  $\forall(u,m)\in \XTM\cap \{1/K\le m\le K,\ |u|^{(1)}_Q \le K\}$, the image $(\hat{u},\hat{m}) :=\cTK(u,m)$ satisfies $1/K\le \hat{m}\le K, |\hat{u}|^{(1)}_Q \le K$.

Suppose $\hat{m}(x,t)$ attain its minimum at $(\hat{x},\hat{t})\in Q$, then from the definition \eqref{eq:Holder-norm-alpha} and \eqref{eq:Holder-norm-alpha+1} of H\"{o}lder spaces and \eqref{eq:bound-hatm-Qnorm-long}, we have 
\begin{equation*}
\frac{|\hat{m}(\hat{x},\hat{t}) - \hat{m}(\hat{x}, 0)|}{T^{1-\frac{d+2}{2r}}} \le \sup_{(x,t_1)\neq (x,t_2)\in Q} \frac{|\hat{m}(x,t_1) - \hat{m}(x,t_2)|}{|t_1-t_2|^{1-\frac{d+2}{2r}}} 
\le |\hat{m}|^{(2-\frac{d+2}{r})}_Q
\le C_2'.
\end{equation*}
Thus, we obtain a lower bound for $\hat{m}$ using \ref{ASS:Init Term Cond} and definition \eqref{eq:K},
\begin{equation}\label{eq:hatm-lower-bound}
\begin{array}{r@{\;}l}
  \min\limits_{(x,t)\in Q}\hat{m}(x,t) = \hat{m}(\hat{x},\hat{t}) & \ge \hat{m}(\hat{x},0) - T^{1-\frac{d+2}{2r}} C_2' = m_0(\hat{x}) - T^{1-\frac{d+2}{2r}} C_2' \\
& \ge \underline{m} - T^{1-\frac{d+2}{2r}} C_2' \ge  2/K - T^{1-\frac{d+2}{2r}} C_2' .
\end{array}
\end{equation}
An analogous argument provides the upper bound for $\hat{m}$,
\begin{equation}\label{eq:hatm-upper-bound}
    \max_{(x,t)\in Q} \hat{m}(x,t) \le |m_0|^{(1)}_{\bbT^d} + T^{1-\frac{d+2}{2r}} C_2' \leq K/2  + T^{1-\frac{d+2}{2r}} C_2' .
\end{equation}
For the estimate on $\hat{u}$, recall from \eqref{eq:bound-hatu-Qnorm} that 
\begin{equation}\label{eq:hatu-bound}
    |\hat{u}|^{(1)}_Q \le |u_T|^{(1)}_{\bbT^d} + T^{\frac12 - \frac{d+2}{2r}} C_3'' \leq K/2  + T^{\frac12 - \frac{d+2}{2r}} C_3''.
\end{equation}
For sufficiently small $T$,  $T^{1-\frac{d+2}{2r}} C_2' \le \min\{1/K, K/2\}$ and $T^{\frac12 - \frac{d+2}{2r}} C_3'' \le K/2$, we conclude from \eqref{eq:hatm-lower-bound}, \eqref{eq:hatm-upper-bound} and \eqref{eq:hatu-bound} that
$1/K\le \hat{m}\le K, |\hat{u}|^{(1)}_Q \le K$. Thus, $\cTK$ is a contraction map on $\XTM \cap \{1/K \le m \le K,\ |u|^{(1)}_Q \le K\}$, and therefore so is $\cT$.

\emph{Step 4: Uniform convergence of policy iteration for inverse MFGs.} To show that the solution sequence $\{b^{(k)} = \cB(q^{(k)}, m^{(k)})\}_{k\ge 0}$ of the policy iteration converges uniformly to a solution $b^*$ of the original inverse problem (i.e., $b^*$ satisfies \eqref{EQ:MFG}). 

To begin with, we look at the sequence $\{(u^{(k)}, m^{(k)})\}_{k\ge 0}$ generated from the policy iteration for inverse MFGs, and show its uniform convergence. Equivalently, we consider the convergence of sequence $\{\cT^{k}(u^{(0)}, m^{(0)})\}_{k\ge 0}$. From the discussions in Step 1-3, $\cT$ is a contraction map on the set $\XTM\cap \{1/K\le m\le K,\ |u|^{(1)}_Q \le K\}$, thus we first need to show the initial iteration $(u^{(0)}, m^{(0)})$ is also in this set.

Given the initialization $q^{(0)}\in C^{1,0}(Q;\bbR^d)$, the initial $(u^{(0)}, m^{(0)})$ satisfies the following equations, obtained by replacing $\Du$ with $q^{(0)}$ in \eqref{EQ:T},
\begin{equation}\label{eq:PDE-policy-q0}
\left\lbrace \begin{array}{@{}l}
\partial_t m^{(0)} - \eps\Delta m^{(0)}- q^{(0)}\cdot \Dhm  - {\rm div}(q^{(0)})m^{(0)} = 0, \; m^{(0)}(\bx,0)=m_0(\bx),\\
-\partial_t u^{(0)}-\eps\Delta u^{(0)} + q^{(0)}\cdot \Dhu - \frac12|q^{(0)}|^2 =
b^{(0)}(\bx)
+F(m^{(0)}), \;
u^{(0)}(\bx,T) = u_T(\bx),
\end{array}\right. 
\end{equation} 
where the initial optimal solution in the inversion is $b^{(0)}(\bx) = -\DataFun(\bx) - \eps\Delta u_T + q^{(0)}(\bx, T)\cdot \Du_T -\frac12 | q^{(0)}(\bx, T)|^2 -F(m^{(0)}(\bx, T))$. Since $q^{(0)}\in C^{1,0}(Q;\bbR^d)$, its regularity and norm estimates provide bounds on the coefficients and source terms in \eqref{eq:PDE-policy-q0}. Following the same ideas from Step 1 by applying the parabolic estimates \Cref{PROP:Parabolic Estimate Strong,PROP:Parabolic Estimate Weak} to the solution $(u^{(0)}, m^{(0)})$ of \eqref{eq:PDE-policy-q0}, we can show for sufficiently small $T$ (depend on $q^{(0)}$ as well), the initial $(u^{(0)}, m^{(0)})\in \XTM$. Following similar arguments as in Step 3, 
we obtain $1/K\le m^{(0)}\le K$ and $|u^{(0)}|^{(1)}_Q \le K$. In conclusion, we have $(u^{(0)}, m^{(0)})\in \XTM\cap \{1/K\le m\le K,\ |u|^{(1)}_Q \le K\}$ when $T$ is sufficiently small.

Therefore, when $T$ is sufficiently small, 
$\cT$ is a contraction on the complete space $\XTM\cap \{1/K\le m\le K,\ |u|^{(1)}_Q \le K\}$ and the initialization $(u^{(0)}, m^{(0)})$ is included in this space. 
From the contraction mapping theorem, $(u^{(k)}, m^{(k)}) = \cT^k(u^{(0)}, m^{(0)})$ converges to a fixed point $(u^*,m^*)$ of $\cT$, i.e.,
\begin{equation}\label{eq:converge-um-fix-point}
\lim_{k\rightarrow\infty} \|u^{(k)} - u^*\|_{W^{2,1}_r(Q)} + |u^{(k)} - u^*|^{(1)}_Q + |m^{(k)} - m^*|^{(1)}_Q = 0.
\end{equation}
Moreover, %
$(u^{(k)}, m^{(k)}), (u^*,m^*) \in\XTM\cap \{1/K\le m\le K,\ |u|^{(1)}_Q \le K\}$.

Next, we consider the convergence of the solution sequence $\{b^{(k)} = \cB(\Du^{(k-1)}, m^{(k)})\}_{k\ge 0}$. 
Recalling the definition of $\cT$ in \eqref{EQ:T}, its fixed point $(u^*, m^*)$ satisfies 
\begin{equation*}
    \left\{\begin{array}{ll}
        \partial_t m^* - \eps\Delta m^* - {\rm div}(m^* \Du^*) = 0\, &\text{in } Q,\\
        -\partial_t u^* - \eps\Delta u^* + \frac12|\Du^*|^2 = \cB(\Du^*, m^*)(\bx) +F(m^*)\, &\text{in } Q,\\
        u^*(\bx,T) = u_T(\bx),\ m^*(\bx,0)=m_0(\bx),\ u^*_t(\bx,T) = \DataFun(\bx)\, &\text{in } \bbT^d.
    \end{array}\right.
\end{equation*}
Therefore, $b^*(\bx) := \cB(\Du^*,m^*)(\bx)$ is a solution to the original inverse problem. It remains only to show that $b^{(k)}(\bx) = \cB(\Du^{(k-1)}, m^{(k)})(\bx)$ converges uniformly to $b^*(\bx)$ on $\bbT^d$. Consider their difference,
\begin{equation}\label{eq:bound-b-diff}
\begin{aligned}
   &\| b^{(k)}(\bx) - b^*(\bx)\|_{L^\infty(\bbT^d)} = \|\cB(\Du^{(k-1)}, m^{(k)})(\bx) - \cB(\Du^*, m^*)(\bx)\|_{L^\infty(\bbT^d)}\\
   & = \|\cB(\psi(\Du^{(k-1)}), \varphi(m^{(k)}))(\bx) - \cB(\psi(\Du^{*}), \varphi(m^*))(\bx)\|_{L^\infty(\bbT^d)}\\
  & \leq C_0(|u^{(k)}-u^*|^{(1)}_Q+ |m^{(k)}-m^*|^{(1)}_Q).
\end{aligned}
\end{equation}
The last line uses $(u^{(k)}, m^{(k)}), (u^*,m^*)\in\XTM\cap \{1/K\le m\le K,\ |u|^{(1)}_Q \le K\}$ and  \eqref{eq:Lipschitz-Tk}.
Taking the limit $k\to\infty$, from \eqref{eq:converge-um-fix-point}, we obtain \eqref{eq:converge-b-fix-point}.
\end{proof}

The proof is mainly based on the parabolic estimates in \Cref{PROP:Parabolic Estimate Strong,PROP:Parabolic Estimate Weak}, by bounding the coefficients and source term of the parabolic equations (such as \eqref{EQ:T}) in the definition of the fixed-point iteration operator. With the discussion \eqref{eq:q-DpHDu} and \eqref{eq:B-simplify}, the theorem and proof can be extended to a general Hamiltonian by adding regularity assumptions on $H(\bp)$.
\edit{Since only the Lipschitz continuity and the boundedness of $F$ are used in the proof, it can be naturally generalized to nonlocal and monotone $F$ with similar smoothness and boundedness assumptions. In addition, the MFG setup can be changed to impose the terminal condition $u(x,T) = F_T(m(x,T))$. In this case, the algorithm and the proof still hold, with additional regularity assumptions required on $F_T$ instead of $u_T$.}

For the cases with different types of data $\DataFun(\bx)$, if we assume the existence of the solution to the linear inverse problem in step (ii), denoted again using the notation $b^{(k)} = \cB(q^{(k)}, m^{(k)})(\bx)$. Additionally, with regularity assumptions on $\cB$, we can derive a similar convergence theorem. However, the existence remains a nontrivial question even for the case (i) with $u(\bx, 0)$ data, which is an inverse parabolic source problem with the final overdetermination \cite{Isakov-Book06,Isakov-Book90,Isakov-CPAM91}.

\subsection{Linear rate of convergence\label{sec:converge-linear}}

We further establish linear convergence of $b^{(k)}$.

\begin{theorem} \label{thm:linear-conv}
	 Under the assumptions \ref{ASS:Init Term Cond}, \ref{ASS:F}, \ref{ASS:Data}, \ref{ASS:q0} and the same setting of \Cref{THM:Convergence}, the sequence $\{ b^{(k)}\}$ generated by the policy iteration method for inverse MFGs has an R-linear rate of convergence, i.e., 
	 there exists constants $\hat{T}$, $C$ and $0<\Gamma<1$ such that for all $T\in (0,\hat{T}]$, the following inequality holds
\begin{equation}\label{eq:converge-R-linear}
    \|b^{(k)} - b^*\|_{L^{\infty}(\bbT^d)} \le C\Gamma^k,\quad \forall k\ge 0.
\end{equation}
\end{theorem}

\noindent\emph{Remark.} For the constants $\hat{T}$, $C$ and $\Gamma$, their dependence is the same as the constant $\bar{T}$ in \Cref{THM:Convergence}, 
namely on the given information of the inverse MFG problems $(\eps, u_T, m_0,F, \DataFun)$ and the initialization $q^{(0)}$.

\begin{proof}[Proof of \Cref{thm:linear-conv}]
	Define $\{(u^{(k)}, m^{(k)}, b^{(k)})\}_{k\ge 0}$ and $(u^*,m^*,b^*)$ the same way as Step 4 in \Cref{THM:Convergence}. Thus, $(u^{(k)}, m^{(k)}) = \cT (u^{(k-1)}, m^{(k-1)})$ and $(u^*, m^*) = \cT (u^*, m^*)$.  For sufficiently small $T$,  we have $(u^{(k)}, m^{(k)}), (u^{(k-1)}, m^{(k-1)}), (u^*,m^*)\in\XTM\cap \{1/K\le m\le K,\ |u|^{(1)}_Q \le K\}$. Moreover, $\cT = \cTK$ is a contraction map in this space from Step 3. Therefore, 
	by applying the contraction argument \eqref{eq:def-contraction} to $(u_1, m_1) = (u^{(k-1)}, m^{(k-1)})$ and $(u_2, m_2) = (u^{*},  m^{*})$,
	the following recursive relation holds, $\forall k\geq 1$,
	\begin{equation}\label{eq:linear-conv-contraction}
	\begin{aligned}
	&\|u^{(k)} - u^*\|_{W^{2,1}_r(Q)} + |u^{(k)} - u^*|^{(1)}_Q + |m^{(k)} - m^*|^{(1)}_Q\\
	&\le \Gamma(\|u^{(k-1)} - u^*\|_{W^{2,1}_r(Q)} + |u^{(k-1)} - u^*|^{(1)}_Q + |m^{(k-1)} - m^*|^{(1)}_Q),
	\end{aligned}
	\end{equation}
	where the constant $\Gamma\in (0,1)$. Thus, applying \eqref{eq:linear-conv-contraction} for all $k\geq 1$ provides
	\begin{equation}\label{eq:converge-um-R-linear}
\begin{aligned}
&\|u^{(k)} - u^*\|_{W^{2,1}_r(Q)} + |u^{(k)} - u^*|^{(1)}_Q + |m^{(k)} - m^*|^{(1)}_Q\\
&\le\Gamma^{k} (\|u^{(0)} - u^*\|_{W^{2,1}_r(Q)} + |u^{(0)} - u^*|^{(1)}_Q + |m^{(0)} - m^*|^{(1)}_Q) \leq C' \Gamma^{k}.
\end{aligned}
\end{equation}
Combined with \eqref{eq:bound-b-diff}, %
$
\| b^{(k)} - b^*\|_{L^\infty(\bbT^d)}  \leq C_0(|u^{(k)}-u^*|^{(1)}_Q+|m^{(k)}-m^*|^{(1)}_Q) \leq C \Gamma^{k}.
$
\end{proof}

R-linear convergence \cite{NoWr-Book06} means the differences between $b^{(k)}$ and $b^*$ are bounded by a sequence $\{C\Gamma^k\}_{k\geq0}$ with a linear rate of convergence. This demonstrates the policy iteration $\{ b^{(k)}\}$ converges exponentially fast to a solution $b^*$ of the original inverse problem. 
Moreover, this R-linear convergence result can be strengthened to standard linear convergence for ${b^{(k)}}$ (i.e., $|b^{(k+1)} - b^*|_{L^{\infty}(\bbT^d)} \le \Gamma |b^{(k)} - b^*|_{L^{\infty}(\bbT^d)}$), under additional regularity assumptions on $F$ and using the relation $|b^{(k)}(\bx) - b^*(\bx)| = |F(m^{(k)}(\bx, T)) - F(m^*(\bx, T))|$ from \eqref{eq:B-simplify}.

\section{Numerical experiments}

In this section, we illustrate the proposed policy iteration method for inverse MFGs using both one-dimensional and two-dimensional examples, and compare its performance with the direct least-squares method.

In the following examples, PDEs are discretized using 
uniform grids in time and space. We use $I$ to denote the number of grid points in space and $N$ as the number of grid points in time.
Following the choices of \cite{CaCaGo-COCV21},
centered second-order finite differences are used for the discrete Laplacian, and rectangular quadrature rules are applied for the integral terms in space.
The Hamiltonian and the divergence term in the FP equation
are both computed via the Engquist-Osher numerical flux for conservation laws, utilizing the two-sided gradient designed to approximate viscosity solutions. An implicit Euler scheme is applied for the time integration. 
All experiments are implemented using \texttt{Matlab}, linear systems are represented using the sparse matrix format \texttt{spdiags}, and the optimization problems are solved using \texttt{fminunc} solver with the quasi-Newton (BFGS) algorithm, providing gradients computed through adjoint equations (details in \cref{sec:adjoint})
if there is no further discussion.

\subsection{Reconstruction of a one-dimensional obstacle function\label{sec:1d}} 
We first consider a one-dimensional problem with the true obstacle function (shown as the yellow solid line in \Cref{fig:1d-reconstruct}), defined as
\begin{equation}\label{eq:b-1d}
b(x) := 0.1(\sin(2\pi x - \sin(4\pi x)) + \exp(\cos(2\pi x))),\, \text{ for } x\in  \bbT,
\end{equation}
which is a smooth function on $\bbT$. We set the final time $T = 1$, the
diffusion coefficient $\eps = 0.3$, the coupling cost $F (m) = m^2 $, and the Hamiltonian $H(\Du) = \frac{1}{2}|\Du|^2$. The initial condition is $m_0(x) = C\exp(-40(x-0.5)^2)$, where $C$ is a normalizing constant ensuring $\int_{\bbT} m_0(x) dx = 1$, and the final condition is $u_T(x) = -m_0(x)$.
To study the performance of the policy iteration for inverse MFGs,
we test our method for two cases: (i) $u(x,0)$ data: Given the information of the initial solution of value function $u$, we want to reconstruct the obstacle function $b(x)$, (ii) $\partial_t u(x,T)$ data: Given the information of the time derivative of value function near final time, we want to reconstruct the obstacle function $b(x)$. Compared with case (i), this data provides more information since it together with final condition $u_T$ provides an extrapolation of $u(x,t)$ among a small period near time $t=T$, but it is also potentially more sensitive to noise.
The data we use is generated from solving \eqref{EQ:MFG} with the true obstacle function \eqref{eq:b-1d}, and is the direct measurement of these MFG solutions through $\MeasOpr u$.
For noisy data, Gaussian noise is added in every discretization point of $\MeasOpr u$ with its magnitude proportional to the $L^2$-norm of the true $\MeasOpr u$.

\begin{figure}[tb]
	\centering
	\begin{tikzpicture}[every mark/.append style={mark size=2pt}]
	\footnotesize
	\begin{groupplot}[group style = {rows=1, columns=3, horizontal
		sep = 36pt, vertical sep=30pt,}, width =  0.24\textwidth, height =
	0.2\textwidth, 
	scale only axis,  
	cycle list name= cyclelist-reconstruct,
	compat = 1.3,
	xlabel={$x$},
	tick label style={font=\scriptsize},
	]
	
	\nextgroupplot[
	legend style = {nodes=right,  font=\scriptsize,
		/tikz/every even column/.append style={column sep=.1cm},
		legend to name=grouplegend-1d-reconstruct},
	legend columns=4,
	title = {$b$},
	]
	\addplot  table[x =x,y=V] {\data/Solve_U0_noise0_gamma0.txt};
	\addlegendentry{Truth }
	\addplot  table[x =x,y=VP] {\data/Solve_U0_noise0_gamma0.txt};
	\addlegendentry{Reconstruction using $u(x,0)$ data }
	\addplot  table[x =x,y=VP] {\data/Solve_UtN_noise0_gamma0.txt};
	\addlegendentry{Reconstruction using $\partial_t u(x,T)$ data }
	
	\nextgroupplot [
	scaled y ticks = false,
	yticklabel style={/pgf/number format/.cd},
	cycle list shift=1,
	title = {Error},
	]
	\addplot  table[x =x,y=Verr] {\data/Solve_U0_noise0_gamma0.txt};
	\addplot  table[x =x,y=Verr] {\data/Solve_UtN_noise0_gamma0.txt};
	\nextgroupplot [
	ymode=log,
	cycle list shift=1,
	xlabel = {\# iteration $k$},
	title = {$\|b^{(k)} - b^*\|_{L^2}$},
	]
	\addplot table[x =ite,y=err] {\data/Ite_U0_noise0_gamma0.txt};
	\addplot table[x =ite,y=err] {\data/Ite_UtN_noise0_gamma0.txt};
	\end{groupplot}
	\node[black] at ($(group c2r1) + (0pt,0.18\textwidth)$) 
	{\pgfplotslegendfromname{grouplegend-1d-reconstruct}}; 
	\end{tikzpicture}
 \vspace{-10pt}
	\caption{Reconstruction results of policy iteration method for the one-dimensional inverse MFG problem \eqref{eq:b-1d}. Left: reconstructed $b$ for difference cases: (i) using $u(x, 0)$ data (blue dashed line) and (ii) using $\partial_t u(x, T)$ (red dotted line), compared with the true obstacle function $b^*$ (yellow solid line). Middle: the absolute error $|b(x) - b^*(x)|$ for different cases. Right: the decay of the error $\|b^{(k)} -b^*\|_{L^2}$ with respect to the number of iterations $k$ (displayed on a logarithmic scale on y-axis).}
	\label{fig:1d-reconstruct}
\end{figure}

For each iteration in the policy iteration algorithm for inverse MFGs, step (i) is to solve the linear Fokker-Planck equation \eqref{eq:FP-policy} forward in time, which corresponds to solving
$N$ linear systems of size $I\times I$. The step (ii) requires solving a linear inverse problem:
the least-squares in case (i) using $u(x,0)$ data is solved using an iterative method (quasi-Newton algorithm), implemented with the \texttt{fminunc} solver in \texttt{Matlab}
with the gradient evaluated using the adjoint method. Each gradient evaluation requires solving an extra linear PDE \eqref{eq:adj-u0-policy}, which is of the same form as the linear FP \eqref{eq:FP-policy} with a different initial condition. The optimization iterations are terminated when the first-order optimality measurement is less than the tolerance, and the number of optimization iterations is around 100 for the first three policy iterations and decreases to around 5 after about ten policy iterations.  The least-squares in case (ii) using $\partial_t u(x,T)$ data can be evaluated directly using  \eqref{eq:B-LS-HJB-UtN}.
Step (iii) in policy iteration can also be directly evaluated since the chosen Hamiltonian leads to $q^{(k+1)} = \Du^{(k)}$. 
The policy iteration algorithm is initialized at $q^{(0)} = 0$, and terminated when the squared distance between policies at successive iterations is below a given tolerance $\tau$, i.e., 
$\max_{t \in [0,T]} \|q^{(k+1)}(\cdot, t)- q^{(k)}(\cdot, t)\|_{L^2(\bbT^d)} < \tau$.  

In \Cref{fig:1d-reconstruct}, we first study the convergence and reconstruction performance of our proposed policy iteration for inverse MFGs with noiseless data. 
Here, we set the tolerance $\tau=10^{-9}$ for policy iteration and $10^{-10}$ for the optimality tolerance in \texttt{fminunc}. For PDE discretization, we select a number of grid points $I = 50$ in space and $N = 100$ grid points in time. The reconstructed $b$ from the policy iteration method for different cases are plotted in \Cref{fig:1d-reconstruct} (Left), compared with the true obstacle function $b^*$ in yellow. Policy iteration methods successfully reconstruct the bump and well structures of the true obstacle function in both cases.
Their absolute errors $|b(x) -b^*(x)|$ are plotted in the middle of \Cref{fig:1d-reconstruct}. These errors are three orders of magnitude smaller than the true values and are clustered around $x=0.5$, where the initial and final conditions $m_0(x)$ and $u_T(x)$ are peaked. In the right of \Cref{fig:1d-reconstruct}, the error $\| b^{(k)} - b^*\|_{L^2}$ with respect to the number of iterations $k$ are plotted in a logarithmic scale on the y-axis. The policy iteration methods require 20-25 iterations to converge, and their errors decay exponentially with respect to the number of iterations, which is consistent with the (R-)linear convergence discussed in \Cref{thm:linear-conv}.

\begin{figure}[tb]
	\centering
		\begin{tikzpicture}[every mark/.append style={mark size=2pt}]
		\footnotesize
	\begin{groupplot}[group style = {rows=1, columns=2, horizontal
		sep = 70pt, vertical sep=30pt,}, width =  0.27\textwidth, height =
	0.2\textwidth, 
	scale only axis,  
	cycle list name= cyclelist-compare,
	compat = 1.3,
	xlabel={Space discretization},
tick label style={font=\scriptsize},
	]
	
	\nextgroupplot[
	ylabel={Time [sec]},
	legend style = {nodes=right, font=\scriptsize,
		/tikz/every even column/.append style={column sep=1cm},
		 legend to name=grouplegend-1d-compare,
	 },
	legend columns=2, 
	title = {Total time cost},
	]
	\addplot  table[x =I,y=tP] {\data/U0_noise0_gamma0.txt};
	\addlegendentry{Policy Iter., $u(x,0)$ data}
	\addplot  table[x =I,y=tD] {\data/U0_noise0_gamma0.txt};
	\addlegendentry{Direct LS, $u(x,0)$ data}
	\addplot table[x =I,y=tP] {\data/UtN_noise0_gamma0.txt};
	\addlegendentry{Policy Iter., $\partial_t u(x,T)$ data}
	\addplot  table[x =I,y=tD] {\data/UtN_noise0_gamma0.txt};
	\addlegendentry{Direct LS, $\partial_t u(x,T)$ data}
	\nextgroupplot [
	scaled ticks=false,tick label style={/pgf/number format/fixed},
	ymax = 0.031,
	ytick={0,0.01,0.02,0.03},
	yticklabels={0,0.01,0.02,0.03},
	ylabel={$\|b - b^*\|_{L^2}/\|b^*\|_{L^2}$},
	title = {Relative reconstruction error},
	]
	\addplot  table[x =I,y=VP] {\data/U0_noise0_gamma0.txt};
	\addplot  table[x =I,y=VD] {\data/U0_noise0_gamma0.txt};
	\addplot table[x =I,y=VP] {\data/UtN_noise0_gamma0.txt};
	\addplot table[x =I,y=VD] {\data/UtN_noise0_gamma0.txt};
	\end{groupplot}
	\node[black] at ($(group c1r1) + (0.22\textwidth,0.19\textwidth)$) 
	{\pgfplotslegendfromname{grouplegend-1d-compare}}; 
	\end{tikzpicture}
 \vspace{-5pt}
	\caption{Comparison of reconstruction time and relative error between policy iteration and direct least-squares method for the one-dimensional inverse MFG problem \eqref{eq:b-1d}. }
	\label{fig:1d-compare}
\end{figure}

In \Cref{fig:1d-compare}, we compare the performance of the policy iteration method for inverse MFG and the direct least-squares method. For the direct LS method, it seeks an obstacle function $b$ to minimize the squared $L^2$-misfit of data; this optimization problem is solved using the quasi-Newton algorithm implemented in the \texttt{fminunc} function, initialized at $b^{(0)}=0$, with gradients provided through solving the corresponding forward-backward coupled adjoint equations (details in \Cref{sec:LS-U0,sec:LS-UtN}) iteratively, similar to the policy iteration for solving the state equation (MFGs). The left figure in \Cref{fig:1d-compare} compares the total computational time for the two methods against the number of grid points in space discretization $I$, while the right figure shows the relative reconstruction error $\|b-b^*\|_{L^2}/\|b^*\|_{L^2}$ of these methods in different cases. Here, we set the tolerance $\tau=10^{-8}$ for policy iteration and fix the number of grid points in time discretization to be $N=100$. Our policy iteration method for inverse MFG is 3-4 times faster than the direct LS method in case (i), while achieving better accuracy. In case (ii), when achieving similar accuracy, our policy iteration method is 10-20 times faster than the direct LS method, and this efficiency difference becomes larger as the number of grid points in the space discretization increases. This comparison demonstrates the superior efficiency and accuracy of our proposed policy iteration method for inverse MFGs, especially its uses in large-scale problems.

The significant reduction in computational time achieved by the policy iteration method for inverse MFGs, compared to the direct LS method, is primarily from the decoupling in the policy iteration method. This decoupling separates the nonlinear optimization with forward-backward coupled PDE constraints into multiple iterations of linear PDE solves and linear inverse problems. The direct LS method relies on the initial choice of $b$ and gradient information. Each objective evaluation in the direct LS method requires solving the MFG state equation (approximately 20 policy iterations, each involving the solutions of 2 PDEs). Additionally, each gradient evaluation requires solving extra adjoint equations (also involving several iterations, each requiring the solutions of 2 PDEs). Thus, each gradient evaluation in the direct LS method entails solving approximately 100 PDEs, which is a considerable computational cost, and must be multiplied by the number of optimization iterations.
In the policy iteration method for inverse MFGs, the optimization step is embedded within step (ii) of the policy iteration. This step is a linear inverse problem, making it significantly easier and cheaper to solve. For case (i), the optimization problem involves a linear PDE constraint and a quadratic objective, meaning each objective evaluation requires solving only one linear PDE, and each gradient evaluation requires solving only one additional PDE. This reduces computational costs compared to the direct LS method, which conversely involves several iterations of PDE system solutions. The computational time for the policy iteration method can be further reduced by using a larger optimization tolerance during the initial policy iterations, thereby saving time in solving the linear inverse problem in step (ii).
For case (ii), the computational savings are even greater because the linear inverse problem is solved in one shot and only requires space discretization of the linear PDE in step (ii).

Note that the time cost and accuracy for case (ii) are much better compared to case (i) for both methods (shown in \Cref{fig:1d-compare}). This difference is mainly because: for noiseless data $\partial_t u(x,T)$ in case (ii) together with the final condition $u(x,T)$ provides more information ($u$ over a short period near $t=T$) compared to the single piece of data $u(x,0)$ in case (i) (which is smoothed out by the diffusion process. \edit{This can be supported by the numerical observation in \Cref{fig:1d-compare-eps} that larger diffusion coefficients lead to larger reconstruction errors in case (i) with $u(x,0)$ as data.})
 Consequently, the inverse problem in case (ii) is easier to solve than in case (i).
 
 \edit{To further study the effect of diffusion coefficients in influencing the convergence rate or the reconstruction error, we compare the reconstruction time and relative error between policy iteration and direct least-squares method for the one-dimensional inverse MFG problem \eqref{eq:b-1d} for different diffusion coefficient $\epsilon$, for fixed period of time $T=1$, in \Cref{fig:1d-compare-eps}. From the numerical results, we observe that larger diffusion coefficients lead to faster convergence of the policy iteration but larger reconstruction errors (especially in case (i)).}
 
 \begin{figure}[tb]
 	\centering
 	\begin{tikzpicture}[every mark/.append style={mark size=2pt}]
 	\footnotesize
 	\begin{groupplot}[group style = {rows=1, columns=2, horizontal
 		sep = 70pt, vertical sep=30pt,}, width =  0.27\textwidth, height =
 	0.2\textwidth, 
 	scale only axis,  
 	cycle list name= cyclelist-compare,
 	compat = 1.3,
 	xlabel={$\epsilon$},
 	tick label style={font=\scriptsize},
 	]
 	
 	\nextgroupplot[
 	ylabel={Time [sec]},
 	legend style = {nodes=right, font=\scriptsize,
 		/tikz/every even column/.append style={column sep=1cm},
 		legend to name=grouplegend-1d-compare-eps,
 	},
 	legend columns=2, 
 	title = {Total time cost},
 	]
 	\addplot  table[x =eps,y=tP] {\data/eps_U0_noise0_gamma0.txt};
 	\addlegendentry{Policy Iter., $u(x,0)$ data}
 	\addplot  table[x =eps,y=tD] {\data/eps_U0_noise0_gamma0.txt};
 	\addlegendentry{Direct LS, $u(x,0)$ data}
 	\addplot table[x =eps,y=tP] {\data/eps_UtN_noise0_gamma0.txt};
 	\addlegendentry{Policy Iter., $\partial_t u(x,T)$ data}
 	\addplot  table[x =eps,y=tD] {\data/eps_UtN_noise0_gamma0.txt};
 	\addlegendentry{Direct LS, $\partial_t u(x,T)$ data}
 	\nextgroupplot [
 	scaled ticks=false,tick label style={/pgf/number format/fixed},
 	ylabel={$\|b - b^*\|_{L^2}/\|b^*\|_{L^2}$},
 	title = {Relative reconstruction error},
 	]
 	\addplot  table[x =eps,y=VP] {\data/eps_U0_noise0_gamma0.txt};
 	\addplot  table[x =eps,y=VD] {\data/eps_U0_noise0_gamma0.txt};
 	\addplot table[x =eps,y=VP] {\data/eps_UtN_noise0_gamma0.txt};
 	\addplot table[x =eps,y=VD] {\data/eps_UtN_noise0_gamma0.txt};
 	\end{groupplot}
 	\node[black] at ($(group c1r1) + (0.22\textwidth,0.19\textwidth)$) 
 	{\pgfplotslegendfromname{grouplegend-1d-compare-eps}}; 
 	\end{tikzpicture}
 	\vspace{-5pt}
 	\caption{\edit{Comparison of reconstruction time and relative error between policy iteration and direct least-squares method for the one-dimensional inverse MFG problem \eqref{eq:b-1d} for different diffusion coefficient $\epsilon$, for fixed period of time $T=1$.} }
 	\label{fig:1d-compare-eps}
 \end{figure}

Furthermore, we address the stability of the proposed policy iteration method for inverse MFG by studying its performance given noisy data. In \Cref{fig:1d-stability}, we show its reconstruction results using $u(x,0)$ data with $1\%$ (in $L^2$-norm) pointwise Gaussian noise. To stabilize the algorithm, we add a Tikhonov-type regularization term $\frac{\gamma}{2}\|\Db\|_{L^2}^2$ to the objective function in the optimization step, where $\gamma$ is the regularization parameter and $\gamma=10^{-6}$ in this case. Other parameters are set the same as in \Cref{fig:1d-reconstruct}.
The blue dashed line in the right figure of \Cref{fig:1d-stability} shows the difference of the reconstructed $b$ compared with the true $b^*$. With only noisy $u(x,0)$ data, the method successfully reconstructs the bump and well structure at the bottom of $b$ (in the area of $x\in[0.4, 0.7]$), but fails to capture the structure near $x=0.2$,  and has a relatively poor reconstruction in area $x>0.9$. Compared with noiseless results in \Cref{fig:1d-reconstruct}, the reconstruction with noisy $u(x,0)$ data also has an error approximately two orders of magnitude larger. Although the reconstructed obstacle function $b$ shows a noticeable difference compared to the true value, the reconstructed initial data $u(x,0)$ (the corresponding solution of MFG given the reconstructed $b$ as input) appears indistinguishable from the true value to the naked eye (shown in the right of \Cref{fig:1d-stability}). This observation highlights the inherent instability of the original inverse MFG problem. \edit{A heuristic interpretation} is that the process is analogous to the inverse heat equation, where the diffusion term smooths out information. As a result, the single piece of data $u(x,0)$ retains only limited information.

\begin{figure}[tbhp]
	\centering
		\begin{tikzpicture}[every mark/.append style={mark size=2pt}]
		\footnotesize
	\begin{groupplot}[group style = {rows=1, columns=3, horizontal
		sep = 36pt, vertical sep=30pt,}, width =  0.24\textwidth, height =
	0.2\textwidth, 
	scale only axis,  
	cycle list name= cyclelist-reconstruct,
	compat = 1.3,
	xlabel={$x$},
	tick label style={font=\scriptsize},
	]
	
	\nextgroupplot[
		legend style = {nodes=right, font=\scriptsize,
		/tikz/every even column/.append style={column sep=.2cm},
		legend to name=grouplegend-u0,
	},
	legend columns=4,
	title = {$b$},
	]
	\addplot  table[x =x,y=V] {\data/Solve_U0_noise0_gamma0.txt};
	\addlegendentry{Truth }
	\addplot  table[x =x,y=VP] {\data/Solve_U0_noise0.01_gamma1e-06.txt};
	\addlegendentry{1\% noisy $u(x,0)$ data }
	\addplot  table[x =x,y=VP] {\data/Solve_U0Uk_noise0.01_gamma9.3277e-05.txt};
	\addlegendentry{1\% noisy $u(x,0)$ and $u(x,0.2)$ data }
		\nextgroupplot [
	cycle list shift=1,
	title = {Error},
	scaled y ticks = false,
	yticklabel style={/pgf/number format/fixed},
	]
	\addplot  table[x =x,y=Verr] {\data/Solve_U0_noise0.01_gamma1e-06.txt};
	\addplot  table[x =x,y=Verr] {\data/Solve_U0Uk_noise0.01_gamma9.3277e-05.txt};
	\nextgroupplot [
	title = {Reconstructed $u(x,0)$},
	]
	\addplot  table[x =x,y=U] {\data/Solve_U0_noise0_gamma0.txt};
	\addplot  table[x =x,y=UP] {\data/Solve_U0_noise0.01_gamma1e-06.txt};
	\addplot  table[x =x,y=UP] {\data/Solve_U0Uk_noise0.01_gamma9.3277e-05.txt};
	\end{groupplot}
	\node[black] at ($(group c2r1) + (0pt,0.18\textwidth)$) 
	{\pgfplotslegendfromname{grouplegend-u0}}; 
	\end{tikzpicture}
 \vspace{-10pt}
	\caption{Reconstruction results of the policy iteration method for the one-dimensional inverse MFG problem \eqref{eq:b-1d} with noisy data. The blue dashed line is for reconstruction using $u(x,0)$ data with $1\%$ noise, and the red dotted line is for reconstruction using extra $u(x,0.2)$ data with $1\%$ noise. Left: reconstructed obstacle function $b$ compared with true $b^*$ (solid yellow line); Middle: error $|b(x)-b^*(x)|$; Right: the corresponding reconstructed $u(x,0)$, from solution of MFG using reconstructed $b$.
		\label{fig:1d-stability}}
\end{figure}

To study the method’s stability with additional data, we also present the reconstruction results using noisy $u(x,0)$ and $u(x,0.2)$ data (both with $1\%$ noise). The reconstruction with the extra data is plotted in \Cref{fig:1d-stability} as the red dotted line and captures the bump and well structure of the true obstacle function more accurately. Although not perfect, it identifies a bent shape near $x=0.2$ that the original reconstruction did not detect and also reconstructs the part near the end with high accuracy. The error is reduced by $30\%$. This demonstrates that adding data improves the stability of the inverse problem, thereby enhancing the stability performance of our proposed policy iteration method for inverse MFGs.

\begin{figure}[tbhp]
	\centering
	\begin{tikzpicture}[every mark/.append style={mark size=2pt}]
	\footnotesize
	\begin{groupplot}[group style = {rows=1, columns=3, horizontal
		sep = 40pt, vertical sep=30pt,}, width =  0.24\textwidth, height =
	0.2\textwidth, 
	scale only axis,  
	cycle list name= cyclelist-reconstruct-longtime,
	compat = 1.3,
	xlabel={$x$},
	tick label style={font=\scriptsize},
	]
	
	\nextgroupplot[
	legend style = {nodes=right,  font=\scriptsize,
		/tikz/every even column/.append style={column sep=.1cm},
		legend to name=grouplegend-1d-reconstruct-longtime},
	legend columns=3,
	title = {$b$},
	]
	\addplot  table[x =x,y=V] {\data/Solve_U0_noise0_gamma0.txt};
	\addlegendentry{Truth }
	\addplot  table[x =x,y=VP] {\data/Solve_U0_noise0_gamma0_T10.txt};
	\addlegendentry{$u(x,0)$ data, $T=10$ }
		\addplot  table[x =x,y=VP] {\data/Solve_U0_noise0_gamma0_T50.txt};
	\addlegendentry{$u(x,0)$ data, $T=50$ }
	\addplot  table[x =x,y=VP] {\data/Solve_UtN_noise0_gamma0_T10.txt};
	\addlegendentry{$\partial_t u(x,T)$ data, $T=10$ }
		\addplot  table[x =x,y=VP] {\data/Solve_UtN_noise0_gamma0_T50.txt};
	\addlegendentry{$\partial_t u(x,T)$ data, $T=50$ }
	
	\nextgroupplot [
	scaled y ticks = false,
	yticklabel style={/pgf/number format/.cd},
	cycle list shift=1,
	title = {Error},
	]

	\addplot  table[x =x,y=Verr] {\data/Solve_U0_noise0_gamma0_T10.txt};
	\addplot  table[x =x,y=Verr] {\data/Solve_U0_noise0_gamma0_T50.txt};
	\addplot  table[x =x,y=Verr] {\data/Solve_UtN_noise0_gamma0_T10.txt};
	\addplot  table[x =x,y=Verr] {\data/Solve_UtN_noise0_gamma0_T50.txt};

	\nextgroupplot [
	ymode=log,
	cycle list shift=1,
	xlabel = {\# iteration $k$},
	title = {$\|b^{(k)} - b^*\|_{L^2}$},
	]
	\addplot table[x =ite,y=err] {\data/Ite_U0_noise0_gamma0_T10.txt};
	\addplot table[x =ite,y=err] {\data/Ite_U0_noise0_gamma0_T50.txt};
	\addplot table[x =ite,y=err] {\data/Ite_UtN_noise0_gamma0_T10.txt};
	\addplot table[x =ite,y=err] {\data/Ite_UtN_noise0_gamma0_T50.txt};

	\end{groupplot}
	\node[black] at ($(group c2r1) + (0pt,0.18\textwidth)$) 
	{\pgfplotslegendfromname{grouplegend-1d-reconstruct-longtime}}; 
	\end{tikzpicture}
	\vspace{-10pt}
	\caption{Reconstruction results of policy iteration method for the one-dimensional inverse MFG problem \eqref{eq:b-1d}, and for longer time periods $T\in\{10, 50\}$. Left: reconstructed $b$ for difference cases and and different terminal time $T$: (i) using $u(x, 0)$ data (blue dashed lines) and (ii) using $\partial_t u(x, T)$ (red dotted lines), compared with the true obstacle function $b^*$ (yellow solid line). Middle: the absolute error $|b(x) - b^*(x)|$ for different cases. Right: the decay of the error $\|b^{(k)} -b^*\|_{L^2}$ with respect to the number of iterations $k$ (displayed on a logarithmic scale on y-axis).}
	\label{fig:1d-reconstruct-longtime}
\end{figure}

\begin{figure}[tbhp]
	\centering
	\begin{tikzpicture}[every mark/.append style={mark size=2pt}]
	\footnotesize
	\begin{groupplot}[group style = {rows=1, columns=2, horizontal
		sep = 70pt, vertical sep=30pt,}, width =  0.27\textwidth, height =
	0.2\textwidth, 
	scale only axis,  
	cycle list name= cyclelist-compare,
	compat = 1.3,
	xlabel={$T$},
	tick label style={font=\scriptsize},
	]
	
	\nextgroupplot[
	ylabel={Time [sec]},
	legend style = {nodes=right, font=\scriptsize,
		/tikz/every even column/.append style={column sep=1cm},
		legend to name=grouplegend-1d-compare-longtime,
	},
	legend columns=2, 
	title = {Total time cost},
	]
	\addplot  table[x =T,y=tP] {\data/T_U0_noise0_gamma0.txt};
	\addlegendentry{Policy Iter., $u(x,0)$ data}
	\addplot  table[x =T,y=tD] {\data/T_U0_noise0_gamma0.txt};
	\addlegendentry{Direct LS, $u(x,0)$ data}
	\addplot table[x =T,y=tP] {\data/T_UtN_noise0_gamma0.txt};
	\addlegendentry{Policy Iter., $\partial_t u(x,T)$ data}
	\addplot  table[x =T,y=tD] {\data/T_UtN_noise0_gamma0.txt};
	\addlegendentry{Direct LS, $\partial_t u(x,T)$ data}
	\nextgroupplot [
	scaled ticks=false,tick label style={/pgf/number format/fixed},
	ylabel={$\|b - b^*\|_{L^2}/\|b^*\|_{L^2}$},
	title = {Relative reconstruction error},
	]
	\addplot  table[x =T,y=VP] {\data/T_U0_noise0_gamma0.txt};
	\addplot  table[x =T,y=VD] {\data/T_U0_noise0_gamma0.txt};
	\addplot table[x =T,y=VP] {\data/T_UtN_noise0_gamma0.txt};
	\addplot table[x =T,y=VD] {\data/T_UtN_noise0_gamma0.txt};
	\end{groupplot}
	\node[black] at ($(group c1r1) + (0.22\textwidth,0.19\textwidth)$) 
	{\pgfplotslegendfromname{grouplegend-1d-compare-longtime}}; 
	\end{tikzpicture}
	\vspace{-5pt}
	\caption{Comparison of reconstruction time and relative error between policy iteration and direct least-squares method for the one-dimensional inverse MFG problem \eqref{eq:b-1d} for longer time horizons $T$. }
	\label{fig:1d-compare-longtime}
\end{figure}

In \Cref{fig:1d-compare-longtime,fig:1d-reconstruct-longtime}, we study the performance of the policy iteration method for inverse MFGs with larger terminal times $T \in \{1, 5, 10, 20, 50, 100\}$. We use the same setup as previous experiments, fixing the spatial discretization at $I=50$ and the time discretization at $N=100T$, keeping the time step constant across different horizons.
\Cref{fig:1d-reconstruct-longtime} shows reconstruction results for longer times $T \in \{10, 50\}$. Compared to \Cref{fig:1d-reconstruct} with $T=1$, we observe that reconstructions using type (i) data $u(\bx, 0)$ have larger errors as $T$ increases, while reconstructions using type (ii) data $\partial_t u(\bx, T)$ remain stable. Similar trends appear in \Cref{fig:1d-compare-longtime}, which compares the performances of policy iteration and the direct LS methods for inverse MFGs varying $T$.
The \edit{heuristic interpretation} is that, for type (i) data, diffusion in the equation gradually dissipates information, leading to higher reconstruction errors over longer horizons. In contrast, type (ii) data provides time-independent information, keeping errors stable as $T$ increases.
Although our theoretical results (\Cref{thm:linear-conv,THM:Convergence}) assume small $T$ for convergence, our numerical tests show that, at least for type (ii) data, the algorithm remains effective for longer time horizons in inverse MFG problems.

\subsection{Reconstruction of a two-dimensional obstacle function}

\begin{figure}[tb]
	\centering
\includegraphics[width=.8\linewidth]{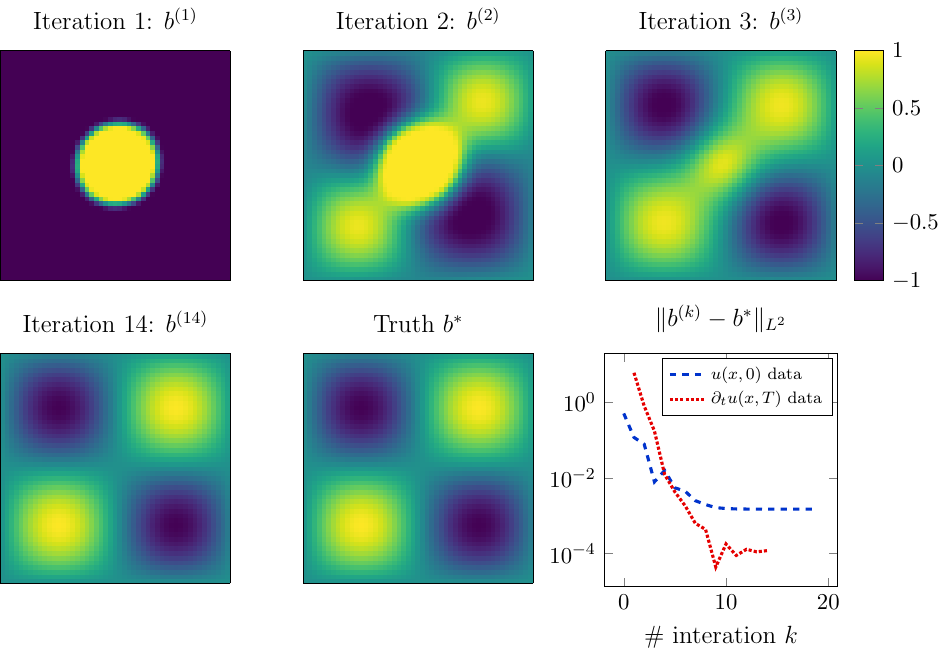}
\vspace{-10pt}
	\caption{Reconstruction results of the policy iteration method for the two-dimensional inverse MFG problem \eqref{eq:b-2d}. The first row shows the reconstructed obstacle function $b^{(k)}$ at different iterations for case (ii): inverse for $b$ using $\partial_t u(\bx, T)$ data. The second row shows the reconstructed obstacle function $b^{(14)}$ at the final iteration and the true obstacle function $b^*$. The last figure is the error $\|b^{(k)}-b^*\|_{L^2}$ versus the number of iterations $k$ for two different cases.
		\label{fig:2d-reconstruct}
	}
\end{figure}

We also test our algorithm in a two-dimensional example to reconstruct the true obstacle function (shown in \Cref{fig:2d-reconstruct}), defined as
\begin{equation}\label{eq:b-2d} 
b(\bx) := \sin(2\pi x_1)\sin(2\pi x_2), \text{ for } \bx=[x_1,x_2]^\top \in  \bbT^2,
\end{equation}
which is a smooth function on $\bbT^2$. Same as the one-dimensional problem, we set the final time $T = 1$, the
diffusion coefficient $\eps = 0.3$, the coupling cost $F (m) = m^2 $, and the Hamiltonian $H(\Du) = \frac{1}{2}|\Du|^2$. The initial condition is $m_0(x) = C\exp(-5[(x_1-0.5)^2+(x_2-0.5)^2])$, where $C$ is a normalizing constant ensuring $\int_{\bbT} m_0(x) dx = 1$, and the final condition is $u_T(x) = -m_0(x)$.
We also study the performance of the policy iteration for inverse MFG,
using two cases: (i) $u(x,0)$ data and (ii) $\partial_t u(x,T)$ data.

In \Cref{fig:2d-reconstruct}, we study the convergence and reconstruction performance of the policy iteration for inverse MFGs in the two-dimensional setting using noiseless data.
Here, we set the tolerance $\tau=10^{-8}$. For the PDE discretization, we select a number of grid points $I = 50\times50=2500$ in space and $N = 100$ grid points in time. For case (ii) using $\partial_t u(\bx, T)$ data, 
the reconstructed $b^{(k)}$ from the policy iteration method at different iterations
are plotted in \Cref{fig:2d-reconstruct}, compared with the true obstacle function $b^*$.
The policy iteration method successfully reconstructs the bump and well structures of the true obstacle function in this example.
The convergence of error $\|b^{(k)}-b^*\|_{L^2}$ versus the number of iterations $k$ is also presented for the two different cases.
Both errors decay exponentially with respect to the number of iterations, which is consistent to the (R-)linear convergence discussed in \Cref{thm:linear-conv}.
The method achieves a three-order of magnitude error decay for case (i).
The error for the $\partial_t u(\bx,T)$ data (case ii) is smaller, due to the same reason we discussed in the one-dimensional case that it provides more data compared with $u(\bx, 0)$ (case i).
We observe that the errors of the policy iteration for inverse MFGs in case (ii) decay rapidly, decreasing by 2-3 orders of magnitude within the first three steps. Combining this with the figures of reconstruction results, we observe that $b^{(1)}$ is dominated by the pattern of the given initial condition $m_0$ and the final condition $u_T$ (bump at the center). By iteration 2, $b^{(2)}$ already detects the patterns in the four corners. By iteration 3, $b^{(3)}$ is already close to the true value, with only a slight difference at the center caused by the initial condition $m_0$ and the final condition $u_T$. Overall, the errors decay by $5-6$ orders of magnitude over $14$ iterations. This demonstrates the accuracy of our proposed policy iteration method for solving inverse MFG in two-dimensional settings.

\begin{figure}[tb]
\centering
	\begin{tikzpicture}[every mark/.append style={mark size=2pt}]
	\footnotesize
\begin{groupplot}[group style = {rows=1, columns=2, horizontal
	sep = 70pt, vertical sep=30pt,}, width =  0.27\textwidth, height =
0.2\textwidth, 
scale only axis,  
cycle list name= cyclelist-compare,
compat = 1.3,
xlabel={Space discretization},
xtick={20, 50, 100},
xticklabels={20$\times$20, 50$\times$50,100$\times$100,},
tick label style={font=\scriptsize},
]

\nextgroupplot[
ylabel={Time [sec]},
legend style = {nodes=right, font=\scriptsize,
	/tikz/every even column/.append style={column sep=1cm},
	legend to name=grouplegend-2d-compare},
legend columns=2, 
title = {Total time cost},
]
\addplot  table[x =Ix,y=tP] {\data/2D_U0_noise0_gamma0.txt};
\addlegendentry{Policy Iter., $u(\bx,0)$ data}
\addplot  table[x =Ix,y=tD] {\data/2D_U0_noise0_gamma0.txt};
\addlegendentry{Direct LS, $u(\bx,0)$ data}
\addplot table[x =Ix,y=tP] {\data/2D_UtN_noise0_gamma0.txt};
\addlegendentry{Policy Iter., $\partial_t u(\bx,T)$ data}
\addplot  table[x =Ix,y=tD] {\data/2D_UtN_noise0_gamma0.txt};
\addlegendentry{Direct LS, $\partial_t u(\bx,T)$ data}
\nextgroupplot [
scaled ticks=false,tick label style={/pgf/number format/fixed},
	ylabel={$\|b - b^*\|_{L^2}/\|b^*\|_{L^2}$},
title = {Relative reconstruction error},
]
\addplot  table[x =Ix,y=VP] {\data/2D_U0_noise0_gamma0.txt};
\addplot  table[x =Ix,y=VD] {\data/2D_U0_noise0_gamma0.txt};
\addplot table[x =Ix,y=VP] {\data/2D_UtN_noise0_gamma0.txt};
\addplot table[x =Ix,y=VD] {\data/2D_UtN_noise0_gamma0.txt};
\end{groupplot}
\node[black] at ($(group c1r1) + (0.22\textwidth,0.19\textwidth)$) 
{\pgfplotslegendfromname{grouplegend-2d-compare}}; 
\end{tikzpicture}
\vspace{-5pt}
	\caption{Comparison of reconstruction time and relative error between policy iteration and direct least-square methods for the two-dimensional inverse MFG problem. }
\label{fig:2d-compare}
\end{figure}

In \Cref{fig:2d-compare}, we compare the performance of the policy iteration method for inverse MFGs and the direct least-squares method. The parameters and conditions are set as previously described.
The left figure in \Cref{fig:2d-compare} compares the total computational time for the two methods against the number of grid points in space discretization $I$, while the right figure shows the relative reconstruction error $\|b-b^*\|_{L^2}/\|b^*\|_{L^2}$ of these methods in different cases. 
Our policy iteration method for inverse MFGs in two-dimensional settings is about five times faster than the direct LS method in case (i) while achieving better accuracy. In case (ii), when achieving similar accuracy, our policy iteration method is 15-25 times faster than the direct LS method, and this efficiency difference becomes larger as the number of grid points in the space discretization increases. 
The main difference from the one-dimensional case (\Cref{fig:1d-compare}) is that the linear system is of squared size due to the discretization in each spatial dimension, which significantly increases the computational burden.
The proposed policy iteration method helps reduce the computational costs, such that even the more challenging case (i) requires less time than the direct LS method for the easier case (ii).
This experiment again demonstrates the superior efficiency and accuracy of our proposed policy iteration method for inverse MFGs, especially its uses in higher-dimensional and large-scale problems.

\section{Conclusion and discussions}

We develop an efficient method based on policy iteration for solving inverse problems in mean-field games. This method (called the policy iteration method for inverse MFGs) reconstructs environmental information (specifically, the obstacle/potential function) from partial observations of value functions. It reformulates the complicated PDE-constrained optimization problem involving strongly coupled MFG systems of nonlinear forward-backward equations into several iterations of linear PDE and linear inverse problem solves, aided by an intermediate variable (the policy). This decoupling accelerates the algorithm through the reduced computational cost of solving linear equations compared to nonlinear coupled PDE systems, as well as their adjoint equations of similar structures needed for gradient computation.
It also simplifies the optimization problem structure by transforming a nonconvex optimization problem with coupled nonlinear PDE constraints into a convex optimization problem with a quadratic objective and linear constraints.
In some special cases, the linear inverse subproblem even has a closed-form solution. We discuss one such case involving derivative data of the value function at the final time. 
From another perspective, the policy iteration method for inverse MFGs simultaneously solves the MFG equations and the original inverse problem, and can be viewed as a fixed-point iteration. We prove the convergence of this approach using the contraction mapping theorem and establish its linear rate of convergence. 
To demonstrate its performance, we compare our approach with the direct least-squares method in reconstructing 1D and 2D obstacle functions. The numerical examples show its superior efficiency (with a 25-times acceleration) and accuracy, and with even better performance on large-scale problems.

There are also opportunities to refine and extend the policy iteration method for inverse MFGs. 
In our numerical experiments, we test the performance of our proposed method for reasonably large time horizons, even though the convergence proof requires small time horizons. Extending the current proof to large $T$ is nontrivial due to complex behaviors that arise, such as the turnpike phenomenon. It is certainly interesting to develop a convergence theory for large $T$ in future work, which will likely require entirely new proof techniques.
In our work, we assume that the data are linear measurements of the value function, primarily because this assumption leads to the linearity of the inverse subproblem in step (ii). However, even without linearity, our approach still has the potential to refine the structure and accelerate the solution of the inverse MFG problem, since the PDE constraint in the inverse subproblem remains linear and propagates in a single direction. Theoretical analysis of this generalized case might still be challenging, as we already face difficulties in the linear case.
Additionally, there is interest in studying potential extensions of this method for other types of data, such as observations of population densities instead of value functions. 
One assumption we make is that $b$ depends only on the state variable $\bx$. While our proposed method can be extended to cases where $b$ also depends on time $t$, this would require significantly more data than the current setup to ensure the uniqueness of the inverse problem.
A more complicated scenario arises when the interaction cost and obstacle are nonseparable (i.e., when the RHS in the first equation of \eqref{EQ:MFG} becomes $\cF(\bx, m)$), and the goal is to recover this operator $\cF$. This leads to a much harder inverse problem, indeed an operator learning problem. In this setting, the policy iteration for MFGs still works, but its application to inversion requires further study.
Moreover, another intriguing extension of the use of the policy iteration method is in the context of Bayesian inverse problems, where the linearity in the inverse subproblem can potentially accelerate Bayesian inversion.

\appendix

\section{Embedding and parabolic estimates}

In this section, we summarize the embedding and parabolic estimates used in the proofs throughout this work, primarily adapted from \cite{LaSoTa-AMO23, CiGiMa-DGA20,BoHaPf-AMO21}. Here, we present only the conclusions of these estimates and omit their proofs, which can be found in the referenced literature.

\begin{lemma}[Lem 2.3 of~\cite{CiGiMa-DGA20}]\label{LEM:Embed 1 1+alpha}
Let $\alpha\in (0,1)$. For any $f\in C^{1+\alpha,(1+\alpha)/2}(Q)$, 
\begin{align*}
    |f|^{(1)}_Q \le |f(\cdot,T)|^{(1)}_{\bbT^d} + T^{\alpha/2}|f|^{(1+\alpha)}_Q,\quad
    |f|^{(1)}_Q \le |f(\cdot,0)|^{(1)}_{\bbT^d} + T^{\alpha/2}|f|^{(1+\alpha)}_Q.
\end{align*}
\end{lemma}

\begin{lemma}[Lem 2.4 of~\cite{CiGiMa-DGA20}]\label{LEM:Embed r 2r} 
Let $r>1$ and $f\in W^{2,1}_{2r}(Q)$. Then 
\begin{equation*}
    \|f\|_{W^{2,1}_r(Q)} \le T^{\frac{1}{2r}}\|f\|_{W^{2,1}_{2r}(Q)}.
\end{equation*}
\end{lemma}

\begin{proposition}[Prop 2.5 of~\cite{CiGiMa-DGA20}]\label{PROP:Embed Holder Sobolev}
Let $f\in W^{2,1}_r(Q)$, and $r$ be such that $r>(d+2)/2$ and $r\neq d+2$. Then 
\begin{equation*}
    |f|^{(2-\frac{d+2}{r})}_Q \le C\left(\|f\|_{W^{2,1}_r(Q)} + \|f(\cdot,0)\|_{W^{2-\frac{2}{r}}_r(\bbT^d)}\right),
\end{equation*}
where $C$ remains bounded for bounded values of $T$. 
\end{proposition}

\begin{lemma}[Lem 2.4 of~\cite{LaSoTa-AMO23}]\label{LEM:Embed Holder Sobolev Cor}
Let $r>d+2$ and $f\in W^{2,1}_r(Q)$. We assume either $f(\bx,0)=0$ or $f(\bx,T)=0$. Then 
\begin{equation*}
    |f|^{(1)}_Q \le CT^{\frac12 - \frac{d+2}{2r}}\|f\|_{W^{2,1}_r(Q)},
\end{equation*}
where $C$ remains bounded for bounded values of $T$.
\end{lemma}

Now consider the linear parabolic problem, 
\begin{equation}\label{EQ:Linear Parabolic}
    \left\{\begin{array}{ll}
        -\partial_t u - \eps\Delta u + b(\bx,t)\cdot \Du + c(\bx,t)u = f(\bx,t)\, &\text{in } Q,\\
        u(\bx,T) = u_T(\bx)\, &\text{in } \bbT^d.
    \end{array}\right.
\end{equation} 
Its solution has the following estimates, bounded by its coefficients and source terms.

\begin{proposition}[Prop 2.7 of~\cite{LaSoTa-AMO23}, Prop 2.6 of~\cite{CiGiMa-DGA20}]\label{PROP:Parabolic Estimate Strong}
Let $r>d+2$ and suppose that $b\in L^{\infty}(Q;\bbR^d)$, $c\in L^{\infty}(Q)$, $f\in L^r(Q)$, and $u_T\in W^{2-\frac{2}{r}}_r(\bbT^d)$. Then the problem \eqref{EQ:Linear Parabolic} admits a unique solution $u\in W^{2,1}_r(Q)$, and it holds that 
\begin{equation}\label{EQ:Parabolic Estimate Strong}
    \|u\|_{W^{2,1}_r(Q)} \le C(\|f\|_{L^r(Q)} + \|u_T\|_{ W^{2-\frac{2}{r}}_r(\bbT^d)}),
\end{equation}
where $C$ depends on the upper bounds for the $L^{\infty}$-norms of the coefficients $b$ and $c$ as well as on $\eps$, $r$, $d$ and $T$, and remains bounded for bounded values of $T$.
\end{proposition}

\begin{proposition}[Prop 2.8 of~\cite{LaSoTa-AMO23}, Thm 4 of~\cite{BoHaPf-AMO21}]\label{PROP:Parabolic Estimate Weak}
Let $r>d+2$ and suppose that $b\in L^{r}(Q;\bbR^d)$, $c\in L^{r}(Q)$, $f\in L^r(Q)$, and $u_T\in W^{2-\frac{2}{r}}_r(\bbT^d)$. Then the problem \eqref{EQ:Linear Parabolic} admits a unique solution $u\in W^{2,1}_r(Q)$, and it holds that 
\begin{equation*}
    \|u\|_{W^{2,1}_r(Q)} \le C,
\end{equation*}
where $C$ depends only on the upper bounds for the $L^{r}$-norms of $b$, $c$, $f$ and $\|u_T\|_{W^{2-\frac{2}{r}}_r(\bbT^d)}$ as well as on $\eps$. Moreover, when $T$ is sufficiently small, the estimate \eqref{EQ:Parabolic Estimate Strong} holds. 
\end{proposition}

\section{Gradient evaluations using adjoint methods \label{sec:adjoint}}

In this section, we provide details of computing gradients using adjoint methods for the PDE-constrained optimization problems discussed in this work.
The adjoint equations presented here are derived at the continuous level and then discretized using the same numerical scheme as the state equations.
 The accuracy of this  optimize-then-discretize approach has been verified using finite difference checks. 
While some discrepancy may exist between the optimize-then-discretize and discretize-then-optimize approaches, our gradient checks show no significant differences (especially for those with linear state equations). 
One can also choose the discretize-then-optimize approach and use the automatic differentiation tools (e.g., \texttt{JAX}) to avoid manual derivation of adjoint equations.
Importantly, the computational cost remains unchanged between the optimize-then-discretize and discretize-then-optimize approaches,  as both approaches involve solving the same number of equations. Therefore, in this paper, we focus on and implement the optimize-then-discretize approach.

\subsection{Adjoint gradient computation for the linear inverse problem in step (ii) of policy iteration with \texorpdfstring{$u(\bx,0)$}{} data}
\label{sec:Adj-Policy-U0}

For the PDE-constrained optimization \eqref{eq:LS-HJB} in step (ii) of policy iteration method for inverse MFGs, we derive its gradient for the case (i) that data $\DataFun(\bx)$ is the observation of value function at the initial time, i.e., $\MeasOpr u:=u(\bx, 0)$. Its Fr\'echet derivative through the adjoint method is:
\begin{equation}\label{eq:grad-u0-policy}
\Phi'(b^{(k)})[\delta b] = \int_{\bbT^d} \left(\int_0^T w(\bx,t)\,dt\right)\delta b(\bx) \,d\bx,
\end{equation}
where $w$ is the solution to the adjoint equation
	\begin{equation}\label{eq:adj-u0-policy}
	\left\lbrace 
\begin{array}{ll}
\partial_t w - \eps\Delta w - {\rm div}(q^{(k)}w) =0\, &\text{in } Q,\\
w(\bx,0) = u^{(k)}(\bx,0) - \DataFun(\bx)\, &\text{in }\bbT^d,
\end{array}\right. 
\end{equation}
where $u^{(k)}$ is the solution to the state equation \eqref{eq:linear-HJB} with the obstacle function on the RHS taking the current value $b^{(k)}$.

\subsection{Adjoint gradient computation for the direct least-squares with \texorpdfstring{$u(\bx,0)$}{} data\label{sec:LS-U0}}

For the PDE-constrained optimization problem \eqref{eq:direct-LS} for the direct least-squares method, we derive its gradient for the case (i) that data $\DataFun(\bx)$ is the observation of value function at the initial time, i.e., $\MeasOpr u:=u(\bx, 0)$. Here, for simplicity, we provide results for $H(\bp) = \frac{1}{2}|\bp|^2$, the similar computation can be extended to a general case of $H$. Its Fr\'echet derivative $\Phi'(b)[\delta b]$ through the adjoint method is the same as in \eqref{eq:grad-u0-policy} with
$(w,v)$ the solution to the adjoint equation \eqref{eq:adj-LS-U0}, which is again a strongly coupled system of forward and backward equations:
\begin{equation}\label{eq:adj-LS-U0}
\left\{\begin{array}{ll}
\partial_t w - \eps\Delta w  - {\rm div}(w \Du) =  {\rm div}(m \Dv)\, &\text{in } Q,\\
- \partial_t v - \eps\Delta v + \Du\cdot \Dv = F'(m)w \, &\text{in } Q,\\
w(\bx,0)= u(\bx,0)-\DataFun(\bx),\, v(\bx,T)=0\, &\text{in } \bbT^d,
\end{array}\right.
\end{equation}
where $(m,u)$ is the solution to the state equation \eqref{EQ:MFG} with the obstacle on the RHS taking values as $b$ in $\Phi'(b)[\delta b]$, and $F'(\cdot)$ is the gradient of $F$.

\subsection{Adjoint gradient computation for the direct least-squares with \texorpdfstring{$\partial_t u(\bx, T)$}{} data\label{sec:LS-UtN}}

For the PDE-constrained optimization problem \eqref{eq:direct-LS} for the direct least-squares method, we derive its gradient for the case (ii) that data $\DataFun(\bx)$ is the observation of derivative information of value function at the final time, i.e., $\MeasOpr u:=\partial_t u(\bx, T)$. Similarly as before, we provide results for $H(\bp) = \frac{1}{2}|\bp|^2$. Its Fr\'echet derivative is:
\begin{center}
$
\Phi'(b)[\delta b] = \int_{\bbT^d} \left[\int_0^T w(\bx,t)\,dt - (\partial_t u(\bx,T) - \DataFun(\bx))\right]\delta b(\bx) \,d\bx,
$
\end{center}
where $(m,u)$ is the solution to the state equation \eqref{EQ:MFG} with the obstacle on the RHS taking values as $b$. $\partial_t u(x,T)$ can be evaluated directly by taking derivative of $u$ with respect to time, or through \eqref{EQ:MFG} at the final time, i.e., $\partial_t u(x,T) = -\eps \Delta u_T + \frac{1}{2}|\Du_T|^2 -b - F(m(x,T))$.
The adjoint variables $(w,v)$ are solution to adjoint equation:
	\begin{equation}\label{eq:adj-LS-UtN}
\left\{\begin{array}{ll}
\partial_t w - \eps\Delta w - {\rm div}(w \Du) = {\rm div}(m \Dv) \, &\text{in } Q,\\
- \partial_t v - \eps\Delta v + \Du\cdot \Dv = F'(m)w\, &\text{in } Q,\\
w(\bx,0) = 0,\, v(\bx,T)= -(\partial_t u(\bx,T)-\DataFun(\bx))F'(m(\bx,T))\, &\text{in }  \bbT^d.
\end{array}\right.
\end{equation}

\section*{Acknowledgments}

This work is partially supported by the National Science Foundation through grants DMS-1937254, and DMS-2309802. The work of S. Tong is additionally supported in part by the the National Science Foundation grant DMS-2529292.

{%
\bibliography{BIB-REN,BIB-PolicyIter}
\bibliographystyle{siamplain}
}

\end{document}